\documentclass[preprint,12pt]{elsarticle}
\usepackage{color}
\usepackage{lineno}
\usepackage{pdflscape}
\usepackage{amsmath}
\usepackage{booktabs}
\usepackage{subfigure}
\usepackage[hyphens]{url}
\usepackage{hyperref}
\usepackage{multirow}
\usepackage{addlines} %A
\usepackage{adjustbox}
\usepackage{siunitx}
\usepackage{url}
\usepackage{natbib} %When changing reference style
\usepackage[utf8]{inputenc}
\usepackage[english]{babel}
\usepackage[dvipsnames]{xcolor}
\usepackage{soul}
\usepackage{adjustbox}
\usepackage{verbatim}

\usepackage{mathtools }

\definecolor{darkblue}{rgb}{0,0.2,0.4}
\usepackage{graphicx}
\usepackage{amssymb}
\usepackage{amsthm}
\usepackage[linesnumbered,ruled,vlined]{algorithm2e}
\usepackage{nomencl}
\usepackage{etoolbox}
\makenomenclature

\renewcommand\nomgroup[1]{%
  \item[\bfseries
  \ifstrequal{#1}{A}{Investment planning related sets}{%
  \ifstrequal{#1}{B}{Operation related sets}{%
  \ifstrequal{#1}{C}{Investment planning related parameters}{%
  \ifstrequal{#1}{D}{Operation related parameters}{%
  \ifstrequal{#1}{E}{Investment planning related variables}{%
  \ifstrequal{#1}{F}{Operation related variables}{
\ifstrequal{#1}{G}{Function}{
  }}}}}}}%
]}

\usepackage{framed} % Framing content
\usepackage{multicol} % Multiple columns environment
\usepackage{eurosym}
\usepackage[utf8]{inputenc}
\usepackage{newunicodechar}

\usepackage[T1]{fontenc}
\usepackage{textcomp}
\usepackage{lmodern}

\usepackage{ltxcmds}

\usepackage{siunitx}
\DeclareSIUnit{\euro}{\texteuro}

\newcommand*{\sieuro}[2][]{\SI[{mode=text,#1}]{#2}{\euro}}

\makeatletter
\newcommand*{\EuroMacro}{}
\protected\def\EuroMacro{%
  \ltx@ifnextchar@nospace\bgroup\sieuro{%
    \ltx@ifnextchar[\sieuro\texteuro
  }%
}
\makeatother

\newcommand{\abbreviations}[1]{%
  \nonumnote{\textit{Abbreviations:\enspace}#1}}

\definecolor{royalpurple}{rgb}{0.58, 0.44, 0.86}

\allowdisplaybreaks
\journal{Energy (Elsevier)}
\journal{Energy Elsevier}
\hypersetup{ colorlinks=true,
linkcolor=darkblue,
filecolor=darkblue,
citecolor=darkblue,
urlcolor=darkblue
}
\begin{document}
\begin{frontmatter}


\title{Modelling and analysis of offshore energy hubs}

% use the tnoteref command within \title for footnotes;
% use the tnotetext command for the associated footnote;
% use the fnref command within \author or \address for footnotes;
% use the fntext command for the associated footnote;
% use the corref command within \author for corresponding author footnotes;
% use the cortext command for the associated footnote;
% \author{Name\corref{cor1}\fnref{label2}}
% \address{Address\fnref{label3}}

% \author[label1,label2]{<author name>}
% \address[label1]{<address>}
% \address[label2]{<address>}
\author[a]{Hongyu Zhang\corref{cor1}}
\author[a]{Asgeir Tomasgard}
\author[b]{Brage Rugstad Knudsen}
\author[b]{Harald G. Svendsen}
\author[a]{Steffen J. Bakker}
\author[c]{Ignacio E. Grossmann}
\cortext[cor1]{Corresponding author's email address: hongyu.zhang@ntnu.no}
\abbreviations{NCS, Norwegian continental shelf; OEH, Offshore energy hub; PFS, Power from shore; Base, The case with only offshore renewables, gas turbines and power from shore; S1, Scenario 1; S2, Scenario 2; MILP, Mixed-integer linear programming.}
\address[a]{Department of Industrial Economics and Technology Management, Norwegian University of Science and Technology, Høgskoleringen 1, 7491, Trondheim, Norway}
\address[b]{SINTEF Energy Research, Kolbjørn Hejes vei 1B, 7491, Trondheim, Norway}
\address[c]{Department of Chemical Engineering, Carnegie Mellon University, 5000 Forbes Avenue, Pittsburgh, PA 15213, USA}

\begin{abstract}
Clean, multi-carrier Offshore Energy Hubs (OEHs) may become pivotal for efficient offshore wind power generation and distribution. In addition, OEHs may provide decarbonised energy supply for maritime transport, oil and gas recovery, and offshore farming while also enabling conversion and temporary storage of liquefied decarbonised energy carriers for export. Here, we investigate the role of OEHs in the transition of the Norwegian continental shelf energy system towards zero-emission energy supply. We develop a mixed-integer linear programming model for investment planning and operational optimisation to achieve decarbonisation at minimum costs. We consider clean technologies, including offshore wind, offshore solar, OEHs and subsea cables. We conduct sensitivity analysis on CO$_2$ tax, CO$_2$ budget and the capacity of power from shore. The results show that (a) a hard carbon cap is necessary for stimulating a zero-emission offshore energy system; (b) offshore wind integration and power from shore can more than halve current emissions, but OEHs with storage are necessary for zero-emission production and (c) at certain CO$_2$ tax levels, the system with OEHs can potentially reduce CO$_2$ emissions by $50\%$ and energy losses by $10\%$, compared to a system with only offshore renewables, gas turbines and power from shore.
%
% Decarbonising the energy supply to installations on the Norwegian Continental Shelf (NCS) is crucial to meet Norway's climate target due to its high share of Norway's total CO$_2$ emission. Switching from gas turbines to renewable energy is one approach and requires investments in offshore renewable technologies and involves multiple energy carriers. System design is essential to achieve a cost-efficient, secure and low emission offshore energy system. Here, we introduce a hybrid energy system operational and investment planning model. This optimisation model aims to provide the best energy mixtures to platforms to achieve decarbonisation at minimum costs. We explore several clean technologies, including offshore floating wind, offshore floating solar, OEH and subsea cables. We conduct sensitivity analysis on three parameters, CO$_2$ tax, CO$_2$ budget and ACC of PFS. In addition, we evaluate the environmental potential of OEH. Case studies are conducted on the North Sea part of the NCS, considering 66 operating fields. We analyse the results using different metrics, including CO$_2$ emissions, energy loss and CO$_2$ marginal price. We find that (a) a hard carbon cap may be necessary for stimulating a zero emission offshore energy system; (b) offshore wind can reduce more than half of the current emission, but OEH may be necessary for zero emission production and (c) the NCS energy system with OEH can potentially reduce up to around 50\% more CO$_2$ emission and 12\% more energy loss than the one without OEH at certain CO$_2$ tax levels.
\end{abstract}
\begin{keyword}
clean offshore energy hub \sep sensitivity analysis \sep deterministic mixed-integer linear programming model 
% MSC codes here, in the form: \MSC code \sep code
% or \MSC[2008] code \sep code (2000 is the default)
\end{keyword}
%Highlights
% Presemt a new vision on creating vir
\end{frontmatter}

% \linenumbers
\section{Introduction}
In 2020, the European Green Deal lifted the EU climate target aiming to reduce greenhouse gas emissions by at least $55\%$ by $2030$ and become the world's first carbon-neutral continent by $2050$ \cite{green_deal}. The North Sea region may play a central role in this with a capacity target of $300$ GW for offshore wind power by $2050$ \cite{eu_wind}. Offshore Energy Hubs (OEHs) and the hub-and-spoke concept offer a transnational and cross-sector solution for better harnessing offshore wind and integration with the rest of the energy system \cite{northsea_wind_pwoerhub}. An energy hub is a physical energy connection point with energy storage where multiple energy carriers can be converted, conditioned \cite{energyhub}. 

Using OEHs to effectively exploit offshore wind power to decarbonise the Norwegian Continental Shelf (NCS) energy system may contribute to meeting Norway's and Europe's climate targets. Norway was the world's third-largest exporter of natural gas in 2019 \cite{bp_norway}. Offshore oil and gas extraction was responsible for $26.6\%$ ($13.3$ Mt CO$_2$ equivalent) of the total Norwegian greenhouse gases in $2020$ \cite{oil_gas_emission}. Norway steps up its climate goal to reduce emissions by $50\%$ -- $55\%$ by $2030$ compared to $1990$ levels \cite{norway_goal}. 

\begin{figure}[t]
    \centering
    \includegraphics[scale=0.7]{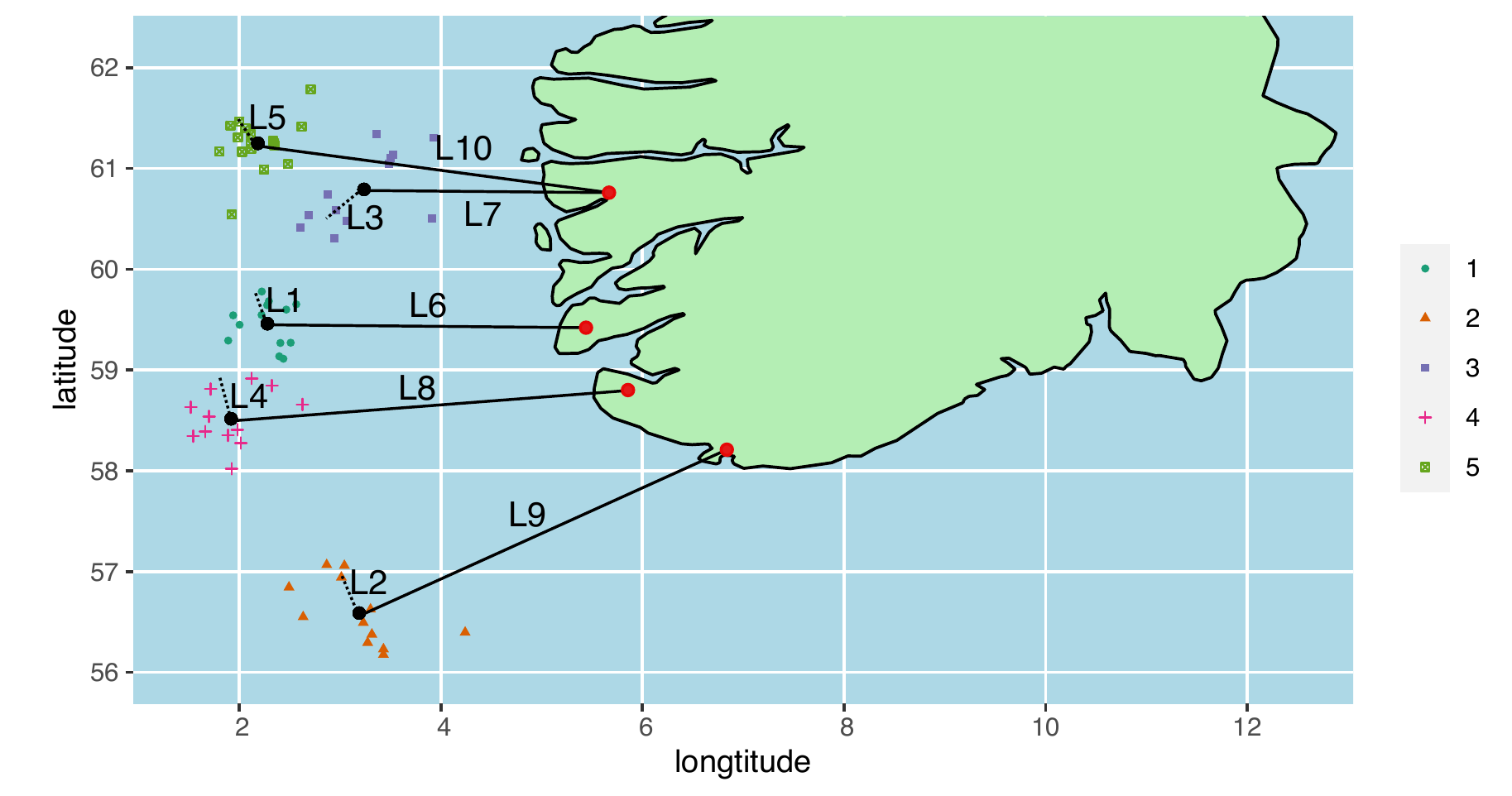}
    \caption{Illustration of the NCS energy system with energy hubs. L1 -- L5 (dotted lines) are representative HVAC cables, while L6 -- L10 (solid lines) are HVDC cables. Black dots represent energy hubs and the red dots represent the onshore buses they connect to. Points with different shapes and colours represent NCS oil and gas fields.}
    \label{system_topology}
\end{figure}

\begin{figure}[!htb]
    \centering
    \includegraphics{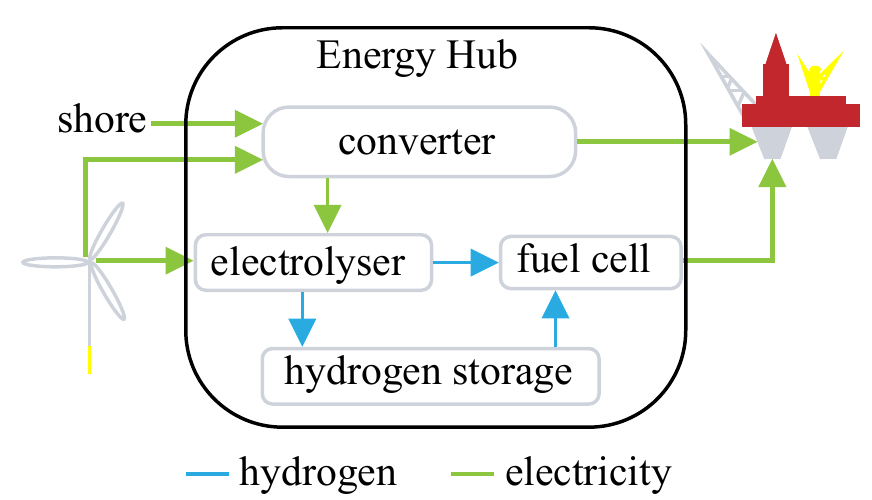}
    \caption{Conceptual illustration of OEHs.}
    \label{offshore_energy_hub}
\end{figure}

This paper investigates the role of clean OEHs in reducing CO$_2$ emissions from the NCS with a focus on NCS oil and gas fields. Figure \ref{system_topology} visualises a potential NCS energy system with OEHs, while Figure \ref{offshore_energy_hub} illustrates the functioning of an OEH. We only consider producing green hydrogen from electrolysis. A likely first use of OEHs may be to provide decarbonised energy supply for offshore oil and gas recovery, maritime cargo transport, and offshore farming \cite{offshore_farming}. In the longer term, the development of OEHs may also facilitate the transition to decarbonised continental energy export. After the lifetimes of oil and gas fields, OEHs may produce and export green hydrogen and connect with a North Sea offshore grid \cite{OEH_greenhydrogen,offshore_grid,northsea_wind_pwoerhub,offshore_book}. 

We develop a multi-carrier Mixed-Integer Linear Programming (MILP) model for investment planning optimisation with a high degree of operational details to cope with the short-term characteristics of the problem. We consider investment in clean technologies, including offshore wind, offshore solar, OEHs and Power From Shore (PFS). Offshore oil and gas platforms involve multiple energy carriers, including electricity, heat, gas, oil and water. More potential energy carriers such as hydrogen and ammonia may become critical during the offshore energy system transition \cite{hydrogen_transition}. The operational model has an hourly time resolution. It is crucial to consider hourly operational decisions in a system with higher penetration of intermittent renewable energy to systematically assess the trade-off between long-term investment decisions and short-term operational decisions \cite{lara_mip}. This is even more critical to be an acceptable alternative to current fossil-based energy provision for installations on the NCS where the security of supply is an absolute requirement. We use a bottom-up approach and model energy demand for devices at oil and gas platforms with hourly time resolution. In this way, we can utilise the potential shifting of loads, such as for water injection, to balance volatile renewable energy better.

We conduct sensitivity analysis on CO$_2$ tax, CO$_2$ budget and the capacity of PFS. Based on the results, we address the value of OEHs in terms of system costs, energy losses and emissions.

The contributions of the paper are: (1) an integrated investment and operational model is developed with the following features: (a) OEHs are modelled for a large-scale offshore energy system and (b) the hourly device-level energy consumption of platforms is modelled, and (2) we analyse the value of OEHs in the North Sea offshore energy system transition towards zero-emission energy supply.

The outline of the paper is as follows: Section \ref{background} introduces the background regarding offshore energy system decarbonisation and offshore oil and gas production. Section \ref{problem_description_modelling strategies_and_assumptions} gives the problem description followed by modelling strategies and assumption. Section \ref{model} presents the model for the case study. Section \ref{results} states the results and analysis of the case study. Section \ref{conclusions} concludes the paper and suggests further research.
\section{Background}
\label{background}
%--------------------------------------------%
This paper tries to find the optimal system design of the NCS energy system for decarbonisation for a future point in time. We consider several decarbonisation solutions, including offshore wind, offshore solar, OEHs and PFS. As our focus is on the NCS, special attention is made to oil and gas recovery. In the following, we present background knowledge on offshore energy system decarbonisation approaches and the production process of offshore platforms.
\subsection{Offshore energy system decarbonisation}
% introduction of four general means of NCS energy system decarbonisation
% =============================================================================
CO$_2$ tax is an important instrument for stimulating offshore energy system decarbonisation. In 2021, the tax is about $55$ €/tonne in Norway with an ambition to increase it to $200$ €/tonne by 2030.
% Norway was one of the first countries to introduce a carbon tax in 1991 \cite{emission_to_air}.
% Another instrument is Greenhouse Gas Emission Trading Act as well as EU Emissions Trading System. 
In addition, the EU Emissions Trading System is a ``cap and trade'' system that also includes the emissions on the NCS \cite{emission_to_air}. Carbon tax and the emissions trading system make a total carbon price of approximately $100$ €/tonne. In this context, oil and gas companies are undertaking considerable investments in decarbonisation solutions to address climate goals, such as PFS and offshore wind. Oil and gas companies on the NCS have set climate targets. For example, Equinor and Vår Energy aim to reduce greenhouse gas emission by $40\%$ by $2030$ and near zero emission by $2050$ \cite{equinor_goal,varenergy_goal}.

% Technologies for decarbonisation exist, and the question is to find the best mixture of such technologies at acceptable costs.
There are four general approaches to reduce offshore CO$_2$ emissions. The first approach is to reduce CO$_2$ emissions by improving reservoir drainage and processing energy efficiency. Water injection and gas injection are common reservoir drainage strategies used on the NCS. 
% In $2019$, $41$ fields on the NCS utilised water injection, and $22$ utilised gas injection \cite{ior}.
% Enhanced oil recovery and water alternating gas are also alternatives. 
Pumping, compression and separation are major processes for handling produced fluids in a processing system. Injection and processing account for more than half of the power consumption at the fields on the NCS. 

Secondly, it is possible to increase the energy efficiency of gas turbines. Due to security of supply requirements, gas turbines usually operate with a margin, which leads to a low efficiency of around $33\%$ \cite{lindegaard}. Adding bottoming cycles to the existing gas turbines can improve their energy efficiency.  However, unlike an onshore energy system, weight and space limitation of an offshore installation restrict extra devices like a bottoming cycle. 

The third approach is to supply zero emission or low emission energy to offshore oil and gas platforms. This includes PFS, switching fuel from natural gas to ammonia or hydrogen and connecting offshore wind farms to platforms.
% ==============================================================================
% review of PFS projects
% ==============================================================================
In the past years, several offshore fields have received PFS via HVDC/HVAC cables  \cite{poweronshore_1,NPD}. 
% are provided with power supply from shore via HVDC/HVAC cables, which includes Ormen Lange, Snøhvit, Troll A, Valhall, Martin Linge, and Goliat \cite{poweronshore_1,NPD}. 
% More fields are planned to receive PFS when the Utsira High Area power grid becomes operational in $2022$ \cite{utsira}.
% However, supplying PFS may have enormous costs. For example, Equinor estimates that the cost to provide power to the Wisting field from shore will be between 380 to 480 € per tonne of carbon dioxide reduced \cite{wisting}, which is highly unprofitable even considering Norway's plan to increase CO$_2$ tax to 200 €/tonne in 2030.

The cost of abating CO$_2$ emissions by taking PFS can vary from less than $100$ to almost $800$ €/tonne \cite{power_from_land}. 
% This is an alternative to building offshore energy solutions, and in the remainder of the paper, we will call it Alternative Abatement Cost (AAC). 
Many abatement projects bringing PFS, are in their planning phase highly unprofitable even considering Norway's plan to increase CO$_2$ tax to $200$ €/tonne in $2030$. Besides, due to the capacity limits of the onshore system, the available power is limited in some cases.

% ==============================================================================
% review of offshore wind research
% ==============================================================================
Offshore wind is another technology to supply clean power to platforms. Equinor's Hywind Tampen project aims to be operational by $2022$ \cite{hywind}. The combination of an offshore platform with a wind farm represents a potential good match for the offshore petroleum sector’s desire for renewable energy with the offshore wind power industry’s desire for an early market \cite{multipleplatforms}. The stability and control issues for an isolated offshore energy system consisting of a wind farm and five platforms were addressed in \cite{multipleplatforms}.
% Following this, they investigated the dynamic performance of wind farm integration to a 5-platform isolated system. 
Integrating large wind turbines into a stand-alone platform is theoretically possible but requires more operational and economic work to prove its feasibility \cite{potential}. In \cite{magnus}, a simulation model that connects a petroleum installation cluster in the Norwegian part of the North Sea to the Norwegian power grid via HVDC cables was presented. They found that local wind power production matching the offshore power demand improves both voltage- and frequency-stability in offshore system.

A framework for determining optimal offshore grid structures for wind power integration and power exchange named Net-Op was presented in \cite{framework}. The authors used a MILP model of the power grid to find optimal grid expansions for connecting offshore wind farms and increasing the transnational power exchange capacity in Europe. They pointed out that a capacity expansion model needs to include a power market model to capture the improved operation of the power system because of the transmission capacity expansion.
% Using the model, they calculated the benefit of an offshore grid from day-ahead spot trade and from the integration of offshore wind power into the system.
An extension of Net-Op that takes into account investment cost, variability of wind/demand/power prices, and the benefit of power trade between countries/price areas was presented in \cite{planningtool}.

The fourth approach is to deploy carbon capture and storage. Storing CO$_2$ in stable underground formations, e.g., old and stable oil reservoirs, has a relatively long history. Since $1996$, nearly one million tonnes of CO$_2$ per year have been separated during the natural gas process from the Sleipner Vest field and stored in the Utsira formation \cite{NPD}.
% The economic and capacity potential of enhanced oil recovery and CO$_2$ storage has been estimated by a techno-economical model consists of a CO$_2$ transportation module and an enhanced oil recovery module that is integrated with an economic model \cite{ccsHOLT}. Simulation of CO$_2$ storage and capacity estimation of four potential CO$_2$ storage hubs is carried out to outline how an expansion in annual CO$_2$ storage capacity offshore the west coast of Norway can be achieved \cite{lothe2019storage}. Improve energy efficiency can reduce CO$_2$ emission up to 30\%. Switch to clean energy sources is necessary to electrify platforms fully. Therefore, in this paper, we focus on the third approach that is supplying clean energy to platforms.

The first two approaches have a limited impact on emission reduction. Whereas the third and fourth approaches can give up to 100\% reduction.

The use of OEHs may allow for better harnessing offshore wind to supply more stable energy to platforms in the short run and export clean energy to the continent in the long run. Using energy hubs for coping with wind power volatility was discussed in \cite{energyhub_wind}. Energy hubs can become sustainable pillars by improving cost-competitiveness of wind power \cite{energyhub_electrofuel}.
% Besides, an energy hub approach was integrated into CGEN model to understand the impact of local energy distribution systems changes on a national scale energy supply system \cite{energyhub_uk}.
In addition, the energy hub concept was also used to increase the energy flexibility in buildings \cite{energyhub_buildings} and electricity markets \cite{energyhub_market}. Techno-economic analysis of offshore energy islands was presented in \cite{northsea_energy}. A separate offshore bidding zone can lead to a more efficient offshore energy system with OEHs \cite{market_arrangement}. An offshore artificial power-to-gas island can produce and transport hydrogen through natural gas pipelines \cite{hydrogen_pipeline}. Energy-hub-based electricity system design for an offshore platform considering CO$_2$ mitigation was presented in \cite{ZHANG20173597}. We refer the readers to \cite{energyhub_review_1, energyhub_review_2} for comprehensive reviews of the research works on energy hubs. 
% In \cite{energyhub_review_2}, the authors showed that electricity and gas networks are the main grids in the reviewed energy hub models, and hydrogen as a clean energy carrier for energy hub can play a crucial role in a sustainable energy supply for various energy sectors.
There is a lack of research on the environmental value of OEHs. Besides research, real-world OEH projects will be carried out \cite{denmark_energyisland, northsea_wind_pwoerhub,poshydon}.

%---------------------------------------
\subsection{Offshore oil and gas fields}
%---------------------------------------
Offshore oil and gas field development and production involve multibillion-dollar investments and profits \cite{babusiaux2007oil}. A large number of studies have been carried out for offshore oil and gas fields production optimisation \cite{production_optimisation} and planning of offshore oil and gas infrastructures \cite{gupta,tarhan}. From the studies, we can see that platforms and fields vary a lot due to, amongst others, geological characteristics, reserves, and remaining lifetimes.
% However, to make large scale investment planning for decarbonising NCS platforms, one may need to simplify the problem based on typical North Sea fields characteristics.
In the following, we present a typical composition and production process of NCS platforms.

\begin{figure}[t]
    \centering
    \includegraphics{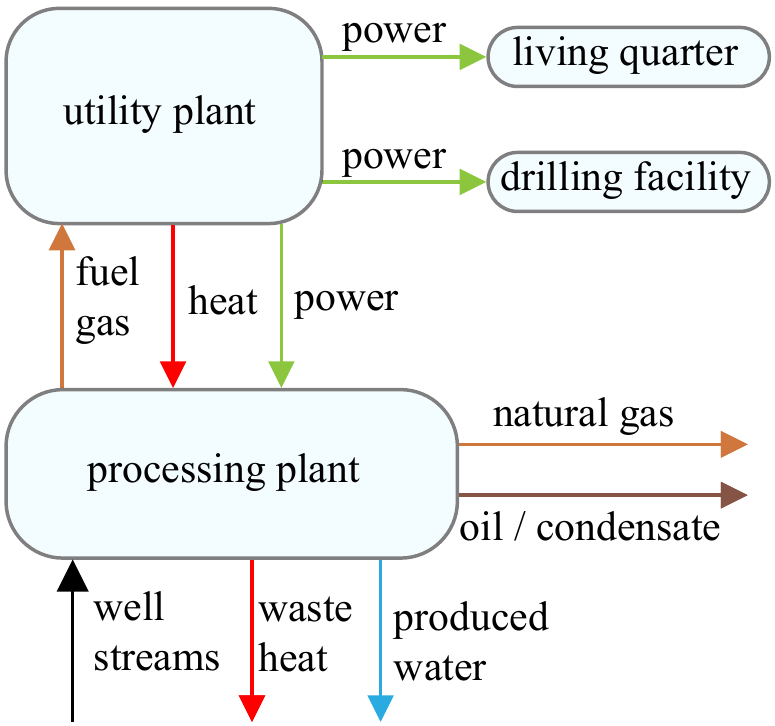}
    \caption{Schematic of a topside structure of a typical North Sea oil and gas platform, adapted from \protect\cite{voldsund}.}
    \label{schematic_platform}
\end{figure}

A North Sea oil field normally consists of topside structures and subsea production systems. A topside structure typically consists of a processing plant, a utility plant, drilling facilities, and a living quarter \cite{voldsund}, see Figure \ref{schematic_platform}. The production plant receives and processes well streams. A visualisation of the production process is presented in  Figure \ref{schematic_production_decarbonisation}. Major energy consumption takes place in the production plants. The energy demand of production plants is conventionally fulfilled by gas turbines located in the utility plant. In $2014$, gas turbines with waste heat recovery units covered approximately $90\%$ of all heat demand for operations on the NCS \cite{WHRU_heat}. 
% The production plant receives and processes well streams, see Figure \ref{schematic_production}. A three-stage separator train separates well streams into produced water, oil, condensate and gas. Typically well streams enter the first separator, which takes out most of the water and gas at arrival conditions. The water content is typically reduced to less than $5\%$ \cite{oil_handbook}. Fuel gas is taken from the first stage separator. The residual mixture of oil, gas and water is usually heated to around $70^{\circ}$C before entering the second stage separator. Produced water is purified and discharged, and in some cases, reinjected into water injection wells to maintain reservoir pressure. Water lift pumps will lift seawater for reinjection if needed. Produced oil is pressurised by pumps and exported via pipelines or shuffle tanks. Produced gas is used as fuel gas, compressed and exported via pipelines or reinjected via dedicated wells for enhanced oil recovery or injected into the same wells for gas lift.

\begin{figure*}[!]
    \centering
    \includegraphics{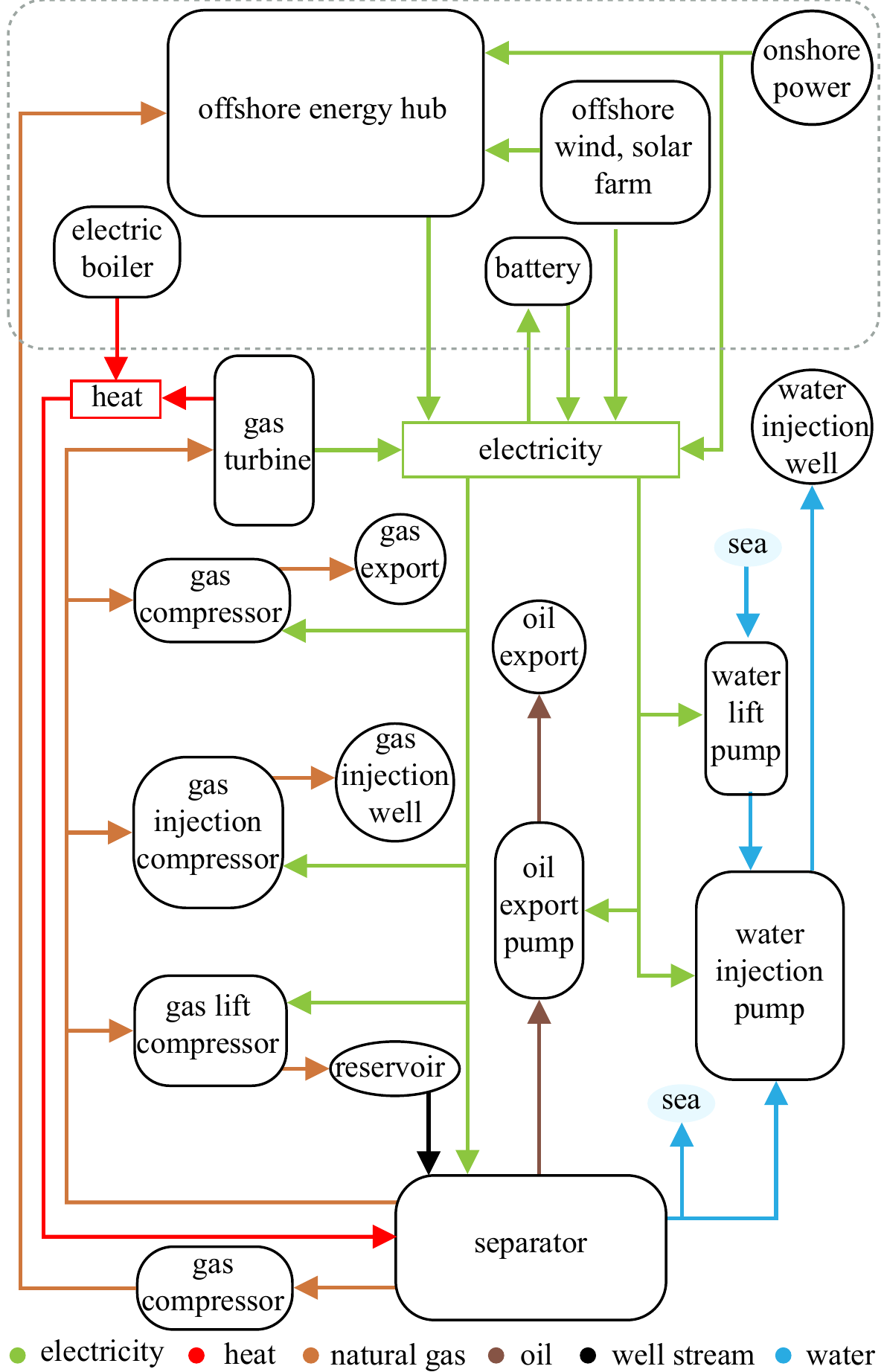}
\end{figure*}
\begin{figure*}[!]
    \caption{Schematic of a potential decarbonised offshore field production process. A three-stage separator train separates well streams into produced water, oil, condensate and gas. Typically the first stage separator takes out most of the water and gas at arrival conditions. Fuel gas is taken from the first stage separator. The residual mix of oil, gas and water is heated before entering the second stage separator. Produced water is purified and discharged, and in some cases, reinjected into water injection wells to maintain reservoir pressure. Water lift pumps will lift seawater for reinjection if needed. Produced oil is pressurised by pumps and exported. Produced gas is used as fuel gas, compressed and exported, reinjected via dedicated wells for enhanced oil recovery or injected into the same wells for gas lift.\\
    The grey dotted box includes the potential processes for decarbonisation. See Figure \protect\ref{offshore_energy_hub} for a visualisation of the processes in an OEH.}
    \label{schematic_production_decarbonisation}
\end{figure*}

\section{Problem description}
\label{problem_description_modelling strategies_and_assumptions}
First, this section introduces the proposed offshore energy system planning problem with OEHs. Then we present the time and geographical structures with the aim of reducing computational time of a potentially big problem. Finally, we give the modelling assumptions.

%\subsection{Problem description}%promote OEH
The problem under consideration aims to make optimal investment and operational decisions for the NCS energy system with OEHs, based on the energy demand captured by the operational model. By solving such a problem, we aim to find out under what conditions OEHs may benefit the system and how OEHs operate with the rest of the system.

To model hourly energy demand, we consider: (a) separators; (b) pumps: water injection pumps, water lift pumps, oil export pumps; (c) compressors: gas injection compressors and gas export compressors. These devices have existing capacities, and no investment will be made in them. Moreover, we assume that device efficiency, flow inlet/outlet pressures and hourly mass flow are given.

For the investments in decarbonisation solutions, we consider: (a) offshore renewable energies (offshore wind and offshore solar); (b) OEHs (electrolysers, hydrogen storage facilities and fuel cells); (c) subsea cables (HVAC, HVDC and offshore and onshore converter stations); (d) electric boilers; (e) platform located batteries. The capital expenditures, fixed operational costs are assumed to be known.

The problem is to determine (a) capacities of decarbonisation technologies and (b) operational strategies that include scheduling of generators, storage and approximate power flow among regions to meet the energy demand with minimum overall investment, operational and environmental costs.
\subsection{Modelling strategies and assumptions}
\label{modelling_strategies_and_assumption}
A multilevel control hierarchy was defined in \cite{foss}, and they argued that the repetitive use of static models could solve many important petroleum production optimisation problems. We develop a multi-period MILP model for an integrated investment planning and operational problem that combines short-term and long-term control hierarchies. We use aggregation, clustering and time sampling \cite{backe} to address the multi-time-scale aspects \cite{kaut} and solve a large-scale instance.
\subsubsection{Time structure of the problem}
\begin{figure}[!thb]
    \centering
    \includegraphics{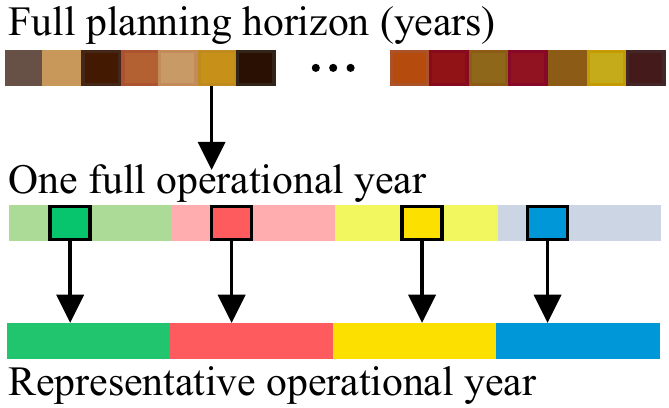}
    \caption{Illustration of combined hierarchies, (adapted from \protect{\cite{time_structure}}).}
    \label{time_structure}
\end{figure}

The investment problem is optimised over a long-term horizon, e.g., a few decades. The operational problem is optimised on an hourly basis based on investment decisions. To combine these two control hierarchies without increasing much the computational time, we select $S$ representative slices, each containing $h$ hours, and scale them up to represent a whole operational year. A visualisation of the time structure is presented in Figure \ref{time_structure}.
\subsubsection{Geographical structure of the problem}
The problem potentially consists of many regions,  and we implement a k-means cluster method based on the locations of fields to reduce the problem size. We have two considerations when deciding the number of clusters. Firstly, we assume the OEH connects the surrounding fields via HVAC cables; thus, only fields with a feasible transmission distance (up to $100$ km) are considered. Secondly, we assume that the cluster centres are the locations for OEHs. We prevent clusters with too few fields. For each cluster, we aggregate the individual fields into one bigger field with a distance to the OEH equal to the average distance of the individual fields and connect fields to OEH in hub-and-spoke form. Currently, we do not consider the interconnection among fields and clusters, resulting in reasonably simple network topology, exemplified in Figure \ref{system_topology}. 

% Gas recompression, gas reinjection, water injection, seawater lift, separation, oil exporting are the processes among others that affect oil and gas production. Compressors, pumps or separators carry out these processes. The processes' energy consumption is affected by flow rate, pressure level and temperature, whose formulations can be quite complex. Therefore, we linearise them in order to conduct a large scale investment-oriented study with considerable operational details.
% We consider hourly time resolution in the operational model. Oil and gas production activity requires a high safety margin, and the operation status may change every second.  Ramping up/down or turning on/off devices can affect energy consumption quite significantly. In addition, wind profile can vary significantly in a very short time. Therefore, the time resolution of the operational model need to be considerable high to guarantee that the investment decisions made can deliver energy under system variations.
\subsubsection{Assumptions}
We assume that each platform is a typical North Sea platform with production processes as presented in Figure \ref{schematic_production_decarbonisation}. The energy consumption of pumps, compressors and separators can be formulated as a function of flow rate, pressure and temperature. For simplicity, we assume the pressure levels and temperatures to take values that are typical on the North Sea, leading to a linear formulation. Kirchhoff voltage law is omitted, and the model is an energy flow model. We assume no mass loss during production.
% meaning that all well streams will be processed with no loss. The produced gas is used for gas turbines power generation, reinjection and exporting. All produced oil is exported by pumps with no loss. The produced water covers the water injection, and the surplus part is discharged overboard. If the produced water is less than water injection, seawater will be lifted for injection by a seawater lift pump. We do not consider the remaining lifetimes of fields.
\section{Mathematical model}
\label{model}
This section presents a deterministic MILP formulation for the multi-carrier energy system investment planning problem with high a degree of operational details. The complete nomenclature of the model can be found in \ref{nomenclature}. The model includes a long-term investment planning horizon and a short-term operational horizon. The integrated investment planning and operational model is partially based upon the LP model developed in \cite{mazzi}. To improve the representation of the fixed capacity independent investment costs, we make use of integer variables. The planner is dedicated to the NCS energy system.
% \subsection{Formulation}
% =============================================================================
% Investment model
% =============================================================================
\subsection{Objective function}
% \subsubsection{Nomenclature}
% \textbf{Sets}\\
\nomenclature[A]{\(\mathcal{P}\)}{set of technologies}
\nomenclature[A]{\(\mathcal{L}\)}{set of transmission branches}
\nomenclature[A]{\(\mathcal{Z}\)}{set of all locations, including platforms, OEH, and onshore buses ($\mathcal{Z}^{P}  \cup \mathcal{Z}^{H}  \cup \mathcal{Z}^{O}$)}
\nomenclature[A]{\(\mathcal{I}\)}{set of operational years}
\nomenclature[A]{\(\mathcal{I}_{0}\)}{set of investment periods}
\nomenclature[A]{\(\mathcal{I}_{i}\)}{set of investment periods $i$ $(i \in \mathcal{I}_{0})$ ancestor to operational year $i$ $(i \in \mathcal{I})$ }
% \textbf{Parameters}\\
\nomenclature[C]{\(C^{PInv}_{pi}\)}{unitary investment cost of technology $p$ at investment period $i$ ($p \in \mathcal{P}, i \in \mathcal{I}_{0}$) [kr/MW, kr/MWh, kr/kg]}
\nomenclature[C]{\(C^{PFix}_{p}\)}{unitary fix operational and maintenance cost of technology $p$ ($p \in \mathcal{P}$) [kr/MW, kr/MWh, kr/kg]}
\nomenclature[C]{\(C^{PFInv}_{pi}\)}{fixed capacity independent investment cost of technology $p$ at investment period $i$ ($p \in \mathcal{P}, i \in \mathcal{I}_{0}$) [kr]}
\nomenclature[C]{\(C^{LInv}_{li}\)}{unitary investment cost of branch $l$ at investment period $i$ ($l \in \mathcal{L}, i \in \mathcal{I}_{0}$) [kr/MW/km]}
\nomenclature[C]{\(C^{LFix}_{l}\)}{unitary fix operational and maintenance cost of branch $l$ ($l \in \mathcal{L}$) [kr/MW/km]}
\nomenclature[C]{\(C^{LFInv}_{li}\)}{fixed capacity independent investment cost of branch $l$ at investment period $i$ ($l \in \mathcal{L}, i \in \mathcal{I}_{0}$) [kr]}
\nomenclature[C]{\(C^{LLInv}_{li}\)}{unitary investment cost of branch $l$ at investment period $i$ ($l \in \mathcal{L}, i \in \mathcal{I}_{0}$) [kr/km]}
\nomenclature[C]{\(X^{PHist}_{pz}\)}{historical capacity of technology $p$ in location $z$ ($p \in \mathcal{P}, z \in \mathcal{Z}$) [MW, MWh, kg]}
\nomenclature[C]{\(X^{PMax}_{pz}\)}{maximum capacity of technology $p$ in location $z$ ($p \in \mathcal{P}, z \in \mathcal{Z}$) [MW, MWh, kg]}
\nomenclature[C]{\(M^{P}_{p}\)}{maximum capacity of a unit of technology $p$ in location $z$ ($p \in \mathcal{P}, z \in \mathcal{Z}$) [MW, MWh, kg]}
\nomenclature[C]{\(\bar{N}_{pz}\)}{maximum number of units of technology $p$ in location $z$ ($p \in \mathcal{P}, z \in \mathcal{Z}$)}
\nomenclature[C]{\(X^{Hist}_{l}\)}{historical capacity of branch $l$ ($l \in \mathcal{L}$) [MW]}
\nomenclature[C]{\(X^{Lgth}_{l}\)}{length of branch $l$ ($l \in \mathcal{L}$) [km]}
\nomenclature[C]{\(X^{LMax}_{l}\)}{maximum capacity of branch $l$ ($l \in \mathcal{L}$) [MW]}
\nomenclature[C]{\(M^{L}_{l}\)}{maximum capacity of a single cable of branch $l$ ($l \in \mathcal{L}$) [MW]}
\nomenclature[C]{\(\bar{N}_{l}\)}{maximum number cables of branch $l$ ($l \in \mathcal{L}$)}
\nomenclature[C]{\(\kappa\)}{scaling effect depending on the number of operation years between investment periods}
% $\pi_{i_{0}}$: probability of investment node $i_{0}$\\
% $\pi_{i}$: probability of operational node $i$\\
% $\mu^{D}_{i}$: relative level of power demand at node $i$\\
% $\mu^{E}_{i}$: yearly CO$_2$ emission limit at node $i$ (t)\\
% $C^{\text{CO}_2}_{i}$: CO$_2$ tax at node $i$ (kr/t)\\
% \textbf{Variables}\\
\nomenclature[E]{\(x_{pzi}^{PAcc}\)}{accumulated capacity of device $p$ in location $z$ in operational year $i$ ($p \in \mathcal{P}, z \in \mathcal{Z}, i \in \mathcal{I}$) [MW, MWh, kg]}
\nomenclature[E]{\(x_{li}^{PAcc}\)}{accumulated capacity of branch $l$ in operational year $i$ ($l \in \mathcal{L}, i \in \mathcal{I}$) [MW]}
\nomenclature[E]{\(x_{pzi}^{PInst}\)}{newly installed capacity of device $p$ in location $z$ at investment period $i_{0}$ ($p \in \mathcal{P}, z \in \mathcal{Z}, i \in \mathcal{I}_{0}$) [MW, MWh, kg]}
\nomenclature[E]{\(x_{li}^{LInst}\)}{newly installed capacity of branch $l$ at investment period $i_{0}$ ($l \in \mathcal{L}, i \in \mathcal{I}_{0}$) [MW]}
\nomenclature[E]{\(\gamma_{pzi}^{PInst}\)}{number of newly installed device $p$ in location $z$ at investment period $i_{0}$ ($p \in \mathcal{P}, z \in \mathcal{Z}, i \in \mathcal{I}_{0}$)}
\nomenclature[E]{\(\gamma_{li}^{LInst}\)}{number of newly installed cables in branch $l$ at investment period $i_{0}$ ($l \in \mathcal{L}, i \in \mathcal{I}_{0}$)}
% $x_{i}$: right hand side of the operational problem\\
% $c_{i}$: cost coefficients of the operational problem\\
% \textbf{Functions}\\
\nomenclature[G]{\(g(\cdot,\cdot)\)}{operational cost in operational year $i$ (kr)}
% The investment planning model is presented in (\ref{mod:investment_planning}).
\begin{equation}
    \min f(\mathbf{x})+\kappa\sum_{i \in \mathcal{I}}g(x_{i}, c_{i})\label{MP_objective}
\end{equation}
The objective function Equation \eqref{MP_objective} is to minimise the total investment ($f(\mathbf{x})$) and operational $(\kappa\sum_{i \in \mathcal{I}}g(x_{i}, c_{i}))$ costs over the planning horizon. The operational cost $g(x_i, c_i)$ is described in Section \ref{sec:operational_model}, where $x_{i}$ and $c_{i}$ are vectors containing capacities and costs information respectively.
\subsection{Investment planning constraints}
\begin{subequations}
\label{mod:investment_planning}
\begin{alignat}{3}
    % &(P1):&&
    % &\min && f(\mathbf{x})+\kappa\sum_{i \in \mathcal{I}}g(x_{i}, c_{i}) \quad& \label{MP_objective}\\
    &&& \mathrlap{f(\mathbf{x})= 
    \sum_{i\in \mathcal{I}_{0}}\left(\sum_{p \in \mathcal{P}} \sum_{z \in \mathcal{Z}}\left(C^{PInv}_{pi}x^{PInst}_{pzi}+C^{PFInv}_{pi}\gamma^{P}_{pzi}\right)\right.}&\notag\\
    &&& \mathrlap{\phantom{f(x) =}\left.+\sum_{l \in \mathcal{L}}\left(C^{LInv}_{li}X^{Lgth}_{l}x^{LInst}_{li}+(C^{LLInv}_{li}X^{Lgth}_{l}+C^{LFinv}_{li})\gamma^{L}_{li}\right) \right) } & \nonumber \\
    &&& \mathrlap{\phantom{f(x) = } +\kappa \sum_{i \in \mathcal{I}}\left(\sum_{p \in \mathcal{P}}\sum_{z \in \mathcal{Z}}C^{PFix}_{pi}x^{PAcc}_{pzi}+\sum_{l \in \mathcal{L}} C^{LFix}_{li}X^{Lgth}_{l}x^{LAcc}_{li}\right) } & \label{investment_cost}\\
    &&&x^{PAcc}_{pzi}=X^{PHist}_{pz}+\sum_{i \in \mathcal{I}_{i}}x^{PInst}_{pzi},& \phantom{abcdefghijkl}  p \in \mathcal{P},  z \in \mathcal{Z},  i \in \mathcal{I} \label{cap_tech}\\\ 
    &&&x^{LAcc}_{li}=X^{LHist}_{l}+\sum_{i \in\mathcal{I}_{i}}x_{li}^{LInst},   & l \in \mathcal{L},  i \in \mathcal{I}\label{cap_line}\\
    &&&0 \leq x^{PInst}_{pzi} \leq \gamma^{P}_{pzi}M^{P}_{p},&  p \in \mathcal{P},  z \in \mathcal{Z},  i \in \mathcal{I}_{0} \label{number_tech}\\
    &&&0 \leq x^{LInst}_{li} \leq \gamma^{L}_{lzi}M^{L}_{l},&  l \in \mathcal{L}, i \in \mathcal{I}_{0} \label{number_line}\\
    &&&0 \leq x^{PAcc}_{pzi} \leq X^{PMax}_{pz},   &  p \in \mathcal{P},  z \in \mathcal{Z},  i \in \mathcal{I} \label{max_tech}\\
    &&&0 \leq x^{LAcc}_{li} \leq X^{LMax}_{l},  & l \in \mathcal{L},  i \in \mathcal{I} \label{max_line}\\
    &&& \sum_{i \in \mathcal{I}_{0}}\gamma^{P}_{pzi} \in \{0, 1, 2, ... , \bar{N}_{pz}\},&p \in \mathcal{P},  z \in \mathcal{Z} \label{tech_number}\\
    &&& \sum_{i \in \mathcal{I}_{0}}\gamma^{L}_{li} \in \{0, 1, 2, ... , \bar{N}_{l}\},&l \in \mathcal{L} \label{line_number}\\
    &&&  \mathrlap{x^{PInst}_{pzi}, x^{LInst}_{li}, x^{PAcc}_{pzi}, x^{LAcc}_{li} \in  \mathbb{R}^{+}_{0},} \label{domain1}\\
    &&& \gamma^{P}_{pzi}, \gamma^{L}_{li} \in \mathbb{Z}^{+}_{0}.
\end{alignat}
\end{subequations}
Equation \eqref{investment_cost} calculates the total capacity dependent investment costs, fixed operating and maintenance costs and fixed capacity independent investment costs. We define $\mathbf{x}$ to be a vector collecting available capacities of all technologies $(\mathcal{P})$ and branches $(\mathcal{L})$ for all operational years $(\mathcal{I})$.  Constraints \eqref{cap_tech} and \eqref{cap_line} represent that the available capacity of a technology $(x^{PAcc}_{pzi})$ or a branch $(x^{LAcc}_{li})$ at an operational year equals to its historical capacity $(X^{PHist}_{pz} \text{ or } X^{LHist}_{l})$ and the sum of newly invested capacities $(x^{PInst}_{pzi} \text{ or } x^{LInst}_{li})$ in its ancestor investment periods $( \mathcal{I}_{i})$. An integer variable $\gamma^{P}_{pzi}$ gives the number of units of technology $p\in \mathcal{P}$, in location $z \in \mathcal{Z}$, in investment period $i \in \mathcal{I}_0$. We define $\gamma^{L}_{li}$ to be the number of cables of branch $l \in \mathcal{L}$. Parameters $M^{P}_{p}$ and $M^{L}_{l}$ state the maximum capacity of a technology unit and a cable, respectively. Parameters $\bar{N}_{pzi}$ and $\bar{N}_{li}$ give the maximum number of units that can be installed of technologies and cables, respectively. 
% Equation \eqref{x} collects all the capacities of all devices and will be passed to the operational model. Equation \eqref{cost} gives the values of parameters, CO$_2$ tax in this case.
% \clearpage
% Different regions have different configurations.
% We make investments in transmission lines, OEH and platform located electric boilers and batteries. HVAC and HVDC lines are distinguished by their costs.
\subsection{Generic operational constraints}
\label{sec:operational_model}
% \subsubsection{Nomenclature}
%------------------------------------------------------------------------%
%---------------------------------- sets --------------------------------%
%------------------------------------------------------------------------%
% \textbf{Sets}
% \begin{itemize}
%     \item General
% \end{itemize}
\nomenclature[B]{\(\mathcal{S}^{T}\)}{set of time slices}
\nomenclature[B]{\(\mathcal{T}\)}{set of hours in all time slices}
% \begin{itemize}
%     \item Electricity
% \end{itemize}
\nomenclature[B]{\(\mathcal{L}\)}{set of subsea cables}
\nomenclature[B]{\(\mathcal{G}\)}{set of gas turbines}
\nomenclature[B]{\(\mathcal{R}\)}{set of renewable units (offshore wind, offshore solar)}
\nomenclature[B]{\(\mathcal{S}^{E}\)}{set of electricity storage}
% \begin{itemize}
%     \item Heat
% \end{itemize}
\nomenclature[B]{\(\mathcal{D}^{Sep}\)}{set of separators}
\nomenclature[B]{\(\mathcal{B}^{E}\)}{set of electric boilers}
% \begin{itemize}
%     \item Natural gas
% \end{itemize}
\nomenclature[B]{\(\mathcal{C}^{Exp}\)}{set of natural gas export compressors}
\nomenclature[B]{\(\mathcal{C}^{Inj}\)}{set of natural gas injection compressors}
% \begin{itemize}
%     \item Oil
% \end{itemize}
\nomenclature[B]{\(\mathcal{P}^{O}\)}{set of oil export pumps}
% \begin{itemize}
%     \item Water
% \end{itemize}
\nomenclature[B]{\(\mathcal{P}^{WI}\)}{set of water injection pumps}
\nomenclature[B]{\(\mathcal{P}^{WL}\)}{set of seawater lift pumps}
% \begin{itemize}
%     \item Hydrogen
% \end{itemize}
\nomenclature[B]{\(\mathcal{E}\)}{set of electrolysers}
\nomenclature[B]{\(\mathcal{F}\)}{set of fuel cells}
\nomenclature[B]{\(\mathcal{S}^{Hy}\)}{set of hydrogen storage facilities}
%------------------------------------------------------------------------%
%------------------------------ parameters -----------------------------%
%------------------------------------------------------------------------%
% \textbf{Parameters}
% \begin{itemize}
%     \item General
% \end{itemize}
\nomenclature[D]{\(C^{\text{CO}_{2}}\)}{CO$_2$ tax (kr/t)}
\nomenclature[D]{\(\mu^{E}\)}{yearly CO$_2$ emission limit (t)}
% $\pi_{i}${probability associated with decision node $i$}
\nomenclature[D]{\(H_t\)}{number of hour(s) in one operational period $t$}
\nomenclature[D]{\(W^{S}_{s}\)}{weight of time slice $s$}
% \begin{itemize}
%     \item Electricity
% \end{itemize}
\nomenclature[D]{\(P_{gz}^{accG}\)}{accumulated capacity of gas turbine $g$ on platform $z$ ($g \in \mathcal{G}, z \in \mathcal{Z}^{P}$) [MW]}
\nomenclature[D]{\(G^{GR}_{g}\)}{ramping limit of gas turbine $g$ ($g \in \mathcal{G}$) [MW/MW]}
\nomenclature[D]{\(R_{rt}^{R}\)}{renewable generation profile of renewable unit $r$ in period $t$ ($r \in \mathcal{R}, t \in \mathcal{T}$)}
\nomenclature[D]{\(\eta_{g}^G\)}{efficiency of gas turbine $g$ ($g \in \mathcal{G}$)}
\nomenclature[D]{\(\eta_{l}^L\)}{efficiency of transmission line $l$ ($l \in \mathcal{L}$)}
\nomenclature[D]{\(\eta_{s}^{SE}\)}{efficiency of electricity store $s$ ($s \in \mathcal{S}^{E}$)}
\nomenclature[D]{\(H_{s}^{SE}\)}{power ratio of electricity store $s$ ($s \in \mathcal{S}^{E}$) [MW/MWh]}
\nomenclature[D]{\(A^{E}_{zl}\)}{elements of bus/line incidence matrix $A^{E}$ ($z \in \mathcal{Z}, l \in \mathcal{L}$)}
\nomenclature[D]{\(\kappa^{PWI}_{p}\)}{electricity demand as fraction of amount of water injected by pump $p$ ($p \in \mathcal{P}^{WI}$) [MW/kg]}
\nomenclature[D]{\(\kappa^{PWL}_{p}\)}{electricity demand as fraction of amount of water lifted by pump $p$ ($p \in \mathcal{P}^{WL}$) [MW/kg]}
\nomenclature[D]{\(\kappa^{PO}_{p}\)}{electricity demand as fraction of amount of oil pumped by pump $p$ ($p \in \mathcal{P}^{O}$) [MW/kg]}
\nomenclature[D]{\(E_{g}^{Fuel}\)}{emission of CO$_2$ of gas turbine $g$ burning fuel ($g \in \mathcal{G}$) [t/MWh]}
\nomenclature[D]{\(C_g^{Fuel}\)}{fuel cost of gas turbine $g$ burning fuel with energy content 1 MWh ($g \in \mathcal{G}$) [kr/MWh]}
\nomenclature[D]{\(C_g^G\)}{variable operational cost of generating 1 MW power from gas turbine $g$ ($g \in \mathcal{G}$) [kr/MW]}
\nomenclature[D]{\(C^{LShed}\)}{load shed penalty cost [kr/MW]}
\nomenclature[D]{\(C^{GShed}\)}{generation shed cost [kr/MW]}
\nomenclature[D]{\(C^{ZO}_{zt}\)}{electricity price in onshore bus $z$ in period $t$ ($z \in \mathcal{Z}^{O}, t \in \mathcal{T} $) [kr/MW]}
\nomenclature[D]{\(\sigma^{Res}_{z}\)}{spinning reserve factor on platform $z$ ($z \in \mathcal{Z}^{P}$)}
% \begin{itemize}
%     \item Heat
% \end{itemize}
\nomenclature[D]{\(P^{accSEP}_{dz}\)}{capacity of separator $d^{sep}$ on platform $z$ ($d \in \mathcal{D}^{Sep}, z \in \mathcal{Z}^{P}$) [MW]}
\nomenclature[D]{\(C^{HLShed}\)}{heat load shed penalty cost [kr/MW]}
\nomenclature[D]{\(\eta^{Gh}_{g}\)}{heat recovery efficiency of gas turbine $g$ ($g \in \mathcal{G}$)}
\nomenclature[D]{\(\eta^{EB}_{b}\)}{efficiency of electric boiler ($b \in \mathcal{B}^{E}$)}
\nomenclature[D]{\(\rho^{SEP}_{d}\)}{heat demand as fraction of amount of oil heated during separation ($d \in D^{Sep}$) [MW/kg]}
% \begin{itemize}
%     \item Natural gas
% \end{itemize}
\nomenclature[D]{\(V_{zt}^{NGE}\)}{natural gas exporting level on platform $z$ in period $t$ ($z \in \mathcal{Z}^{P}, t \in \mathcal{T}$) [kg]}
\nomenclature[D]{\(V_{zt}^{NGI}\)}{natural gas injection level on platform $z$ in period $t$ ($z \in \mathcal{Z}^{P}, t \in \mathcal{T}$) [kg]}
\nomenclature[D]{\(P^{accCExp}_{cz}\)}{capacity of compressor $c$ for exporting on platform $z$ ($c \in \mathcal{C}^{Exp}$, $z \in \mathcal{Z}^{P}$) [MW]}
\nomenclature[D]{\(P^{accCInj}_{cz}\)}{capacity of compressor $c$ for injection on platform $z$ ($c \in \mathcal{C}^{Inj}$, $z \in \mathcal{Z}^{P}$) [MW]}
\nomenclature[D]{\(\gamma_{c}\)}{compression ratio of compressor $c$ ($c \in \mathcal{C}^{Inj}  \cup \mathcal{C}^{Exp}$) }
% $\gamma_{c^{exp}}\)}{compression ratio of export compressor $c^{exp}$}
\nomenclature[D]{\(\alpha\)}{polytropic exponent of compressor}
\nomenclature[D]{\(\eta^{C}\)}{efficiency of compressor}
% $\theta^{NG}$: natural gas energy content (MWh/kg)
% \begin{itemize}
%     \item Oil
% \end{itemize}
\nomenclature[D]{\(P_{pz}^{accPO}\)}{capacity of oil pump $p$ on platform $z$ ($p \in \mathcal{P}^{O}$, $z \in \mathcal{Z}^{P}$) [MW]}
\nomenclature[D]{\(V^{OD}_{zt}\)}{oil exporting level on platform $z$ in period $t$ ($z \in \mathcal{Z}^{P}$, $t \in \mathcal{T}$) [kg]}
% \begin{itemize}
%     \item Water
% \end{itemize}
\nomenclature[D]{\(P_{pz}^{accWL}\)}{capacity of water lift pump $p$ on platform $z$ ($p \in \mathcal{P}^{WI}$, $z \in \mathcal{Z}^{P}$) [MW]}
\nomenclature[D]{\(P_{pz}^{accWI}\)}{capacity of water injection pump $p$ on platform $z$ ($p \in \mathcal{P}^{WL}$, $z \in \mathcal{Z}^{P}$) [MW]}
\nomenclature[D]{\(V^{WI}_{zt}\)}{water injection level on platform $z$ in period $t$  ($z \in \mathcal{Z}^{P}$, $t \in \mathcal{T}$) [kg]}
\nomenclature[D]{\(V^{WL}_{zt}\)}{water lift level on platform $z$ in period $t$  ($z \in \mathcal{Z}^{P}$, $t \in \mathcal{T}$) [kg]}
\nomenclature[D]{\(V^{WB}_{zt}\)}{produced water on platform $z$ in period $t$  ($z \in \mathcal{Z}^{P}$, $t \in \mathcal{T}$) [kg]}
% \begin{itemize}
%     \item Hydrogen
% \end{itemize}
\nomenclature[D]{\(F^{FR}_{f}\)}{ramping limit of fuel cell $f$ ($f \in \mathcal{F}$) [MW/MW]}
\nomenclature[D]{\(\eta^{ES}\)}{conversion factor of electrolyser to inject hydrogen to the storage facility [MWh/kg]}
\nomenclature[D]{\(\eta^{EF}\)}{conversion factor of electrolyser to inject hydrogen directly to fuel cell [MWh/kg]}
\nomenclature[D]{\(\eta^{F}_{f}\)}{efficiency of fuel cell $f$ ($f \in \mathcal{F}$)}
\nomenclature[D]{\(\theta^{Hy}\)}{hydrogen energy content [MWh/kg]}
%------------------------------------------------------------------------%
%------------------------------ variables -----------------------------%
%------------------------------------------------------------------------%
% \noindent \textbf{Variables}
% \begin{itemize}
%     \item Electricity
% \end{itemize}
% \nonmenclature{p_{rzi}^{accR}}{accumulated capacity of renewable unit $r$ near OEH $z$ at node $i$ ($r \in \mathcal{R}, z \in \mathcal{Z}^{H}$) [MW]}\\
\nomenclature[F]{\(q_{sz}^{accSE}\)}{accumulated storage capacity of electricity store $s$ on platform $z$ ($s \in \mathcal{S}^{E}, z \in \mathcal{Z}^{P}$) [MWh]}
\nomenclature[F]{\(p_{z^{o}}^{accZO}\)}{accumulated capacity of shore bus $z$ ($z \in \mathcal{Z}^{O}$) [MW]}
\nomenclature[F]{\(p_{l}^{accL}\)}{accumulated capacity of line $l$ ($l \in \mathcal{L}$) [MW]}
\nomenclature[F]{\(p_{zt}^{D}\)}{total power demand on platform $z$ in period $t$ ($z \in \mathcal{Z}^{P}$) [MW]}
\nomenclature[F]{\(p_{gzt}^{G}\)}{power output of gas turbine $g$ on platform $z$ in period $t$ ($g \in \mathcal{G}, z \in \mathcal{Z}^{P}, t \in \mathcal{T}$) [MW]}
\nomenclature[F]{\(p_{gzt}^{ResG}\)}{power reserved of gas turbine $g$ for spinning reserve requirement on platform $z^{p}$ in period $t$ ($g \in \mathcal{G}, z \in \mathcal{Z}^{P}, t \in \mathcal{T}$) [MW]}
\nomenclature[F]{\(p_{rzt}^{R}\)}{power output of renewable unit $r$ near energy hub $z$ in period $t$ ($r \in \mathcal{R}, z \in \mathcal{Z}^{H}, t \in \mathcal{T}$) [MW]}
\nomenclature[F]{\(p_{szt}^{SE+(-)}\)}{charge (discharge) power of electricity store $s$ on platform $z$ in period $t$ ($s \in \mathcal{S}^{E}, z \in \mathcal{Z}^{P}, t \in \mathcal{T}$) [MW]}
\nomenclature[F]{\(p_{szt}^{ResSE}\)}{power reserved in electricity store $s$ for spinning reserve requirement on platform $z$ in period $t$  ($s \in \mathcal{S}^{E}, z \in \mathcal{Z}^{P}, t \in \mathcal{T}$) [MW]}
\nomenclature[F]{\(q_{szt}^{SE}\)}{energy level of electricity store $s$ on platform $z$ at the start of period $t$ ($s \in \mathcal{S}^{E}, z \in \mathcal{Z}^{P}, t \in \mathcal{T}$) [MWh]}
\nomenclature[F]{\(p_{zt}^{GShed}\)}{electricity generation shed at $z$ in period $t$ ($z \in \mathcal{Z}, t \in \mathcal{T}$) [MW]}
\nomenclature[F]{\(p_{zt}^{LShed}\)}{electricity load shed on platform $z$ in period $t$ ($z \in \mathcal{Z}^{P}, t \in \mathcal{T}$) [MW]}
\nomenclature[F]{\(p_{lt}^{L}\)}{power flow in line $l$ in period $t$ ($l \in \mathcal{L}, t \in \mathcal{T}$) [MW]}
\nomenclature[F]{\(p_{zt}^{ZO}\)}{power supply from onshore bus $z$ in period $t$ ($z \in \mathcal{Z}^{O}, t \in \mathcal{T}$) [MW]}
% \begin{itemize}
%     \item Heat
% \end{itemize}
\nomenclature[F]{\(p^{accEB}_{bz}\)}{accumulated capacity of electric boiler $b$ on platform $z$ ($b \in \mathcal{B}^{E}, z \in \mathcal{Z}^{P}$) [MW]}
\nomenclature[F]{\(p_{bzt}^{EB}\)}{electricity input to electric boiler $b$ on platform $z$ in period $t$ ($b \in \mathcal{B}^{E}, z \in \mathcal{Z}^{P},t \in \mathcal{T}$) [MW]}
\nomenclature[F]{\(p_{dzt}^{HSEP}\)}{heat demand of separator $d$ on platform $z$ in period $t$ ($d \in \mathcal{D}^{Sep}, z \in \mathcal{Z}^{P}, t \in \mathcal{T}$) [MW/kg]}
\nomenclature[F]{\(p_{zt}^{HGShed}\)}{heat generation shed on platform $z$ in period $t$ ($z \in \mathcal{Z}^{P}, t \in \mathcal{T}$) [MW]}
\nomenclature[F]{\(p_{zt}^{HLShed}\)}{heat load shed on platform $z$ in period $t$ ($z \in \mathcal{Z}^{P}, t \in \mathcal{T}$) [MW]}
% \begin{itemize}
%     \item Natural gas
% \end{itemize}
\nomenclature[F]{\(p_{czt}^{CExp}\)}{power demand of the compressor $c$ on platform $z$ in period $t$ ($c \in \mathcal{C}^{Exp}, z \in \mathcal{Z}^{P}, t \in \mathcal{T}$) [MW]}
\nomenclature[F]{\(p_{czt}^{CInj}\)}{power demand of the compressor $c$ on platform $z$ in period $t$ ($c \in \mathcal{C}^{Inj}, z \in \mathcal{Z}^{P}, t \in \mathcal{T}$) [MW]}
% \begin{itemize}
%     \item Oil
% \end{itemize}
\nomenclature[F]{\(p_{pzt}^{PO}\)}{oil exporting pump power $p$ power consumption on platform $z$ in period $t$ ($p \in \mathcal{P}^{O}, z \in \mathcal{Z}^{P}, t \in \mathcal{T}$) [MW]}
% \begin{itemize}
%     \item Water
% \end{itemize}
\nomenclature[F]{\(p_{pzt}^{PWI}\)}{water injection pump $p$ power consumption on platform $z$ in period $t$ ($p \in \mathcal{P}^{WI}, z \in \mathcal{Z}^{P}, t \in \mathcal{T}$) [MW]}
\nomenclature[F]{\(p_{pzt}^{PWL}\)}{seawater lift pump $p$ power consumption on platform $z$ in period $t$ ($p \in \mathcal{P}^{WL}, z \in \mathcal{Z}^{P}, t \in \mathcal{T}$) [MW]}
% \begin{itemize}
%     \item Hydrogen
% \end{itemize}
\nomenclature[F]{\(p_{fz}^{accF}\)}{accumulated capacity of fuel cell $f$ on platform $z$ ($f \in \mathcal{F}, z \in \mathcal{Z}^{H}, t \in \mathcal{T}$) [MW]}
\nomenclature[F]{\(p_{ez}^{accE}\)}{accumulated capacity of electrolyser $e$ on platform $z$ ($e \in \mathcal{E}, z \in \mathcal{Z}^{H}, t \in \mathcal{T}$) [MW]}
\nomenclature[F]{\(v_{sz}^{accSHy}\)}{accumulated storage capacity of hydrogen store $s$ on platform $z$ ($s \in \mathcal{S}^{Hy}, z \in \mathcal{Z}^{H}, t \in \mathcal{T}$) [kg]}
\nomenclature[F]{\(p^{F}_{fzt}\)}{power output of fuel cell $f$ on OEH $z$ in period $t$ ($f \in \mathcal{F}, z \in \mathcal{Z}^{H}, t \in \mathcal{T}$) [MW]}
\nomenclature[F]{\(p^{E}_{ezt}\)}{power input of electrolyser $e$ on OEH $z$ in period $t$ ($e \in \mathcal{E}, z \in \mathcal{Z}^{H}, t \in \mathcal{T}$) [MW]}
% $v^{EF}_{ez^{h}t}\)}{hydrogen from electrolyser $e$ to fuel cell $f$ on OEH $z^{h}$ in period $t$ [kg]}
\nomenclature[F]{\(v^{SHy+(-)}_{szt}\)}{injection (withdraw) of hydrogen to (from) hydrogen storage $s$ on OEH $z$ in period $t$ ($s \in \mathcal{S}^{Hy}, z \in \mathcal{Z}^{H}, t \in \mathcal{T}$) [kg]}
% $H_{s^{h_{y}}}^{SHy}\)}{ratio of storage level capacity and injection [withdraw] capacity of hydrogen storage $s^{h_{y}}$ [kg/kg]}
% $v^{DHy}_{z^{h}t}\)}{hydrogen demand on platform $z^{h}$ in period $t$ [kg]\\
\nomenclature[F]{\(v^{SHy}_{szt}\)}{storage level of hydrogen storage $s$ on OEH $z$ in period $t$ ($s \in \mathcal{S}^{Hy}, z \in \mathcal{Z}^{H}, t \in \mathcal{T}$) [kg]}
% \subsubsection{Operational model objective function}
% \label{model}
We now present the generic operational constraints in one operational year $i$. Note that we omit index $i$ in the operational model for ease of notation.
% The MILP investment model can be a potentially very difficult problem to solve. Embedding more operational constraints in $(P1)$ makes it more complex. Therefore, we separate the operational model formulation $(P2)$ from $(P1)$ for the possibility of applying decomposition techniques \cite{benders} in future.
We model oil and gas recovery as this is the most likely use in the short to medium term. The operational constraints can be modified for other use, e.g., offshore fish farming, maritime, transport, and others. The complete description of operational constraints is presented in \ref{full_operational_model}.
\begin{subequations}
\label{mod:operational}
\begin{alignat}{3}
    &\mathrlap{g(x,c)=\sum_{s \in \mathcal{S}^{T}}\sum_{t \in \mathcal{T}}\sum_{z \in \mathcal{Z}}W_{s}\left(\sum_{p \in \mathcal{P}}C^{P}_{p}y^{P}_{pzt}+C^{Shed}y^{LShed}_{zt}+ \right. } \notag \\ 
    & \mathrlap{\hspace{9cm} \left. \sum_{g \in \mathcal{G}}C^{\text{CO}_2}\rho^{E}_{g}y^{G}_{gzt} \right)}& \label{SP_objective}\\
    % &\mathrlap{\min\sum_{s \in \mathcal{S}^{T}}\sum_{t \in \mathcal{T}}\sum_{z \in \mathcal{Z}}W_{s}\left(\sum_{p \in \mathcal{P}}C^{P}_{p}y^{P}_{pzt}+C^{Shed}y^{LShed}_{zt}+\sum_{g \in \mathcal{G}}C^{\text{CO}_2}\rho^{E}_{g}y^{G}_{gzt} \right)}& \label{SP}\tag{P2}\\
     &&&0 \leq y^{P}_{pzt} \leq x^{PAcc}_{pz}, & p \in \mathcal{P},  z \in \mathcal{Z}, t \in \mathcal{T}\label{gmod:tech_cap}\\
     &&& -x^{LAcc}_{l} \leq y^{L}_{lt} \leq x^{LAcc}_{l}, & l \in \mathcal{L}, t \in \mathcal{T}\label{gmod:line_cap}\\
    &&&\mathrlap{\sum_{g \in \mathcal{G}}y^{G}_{gt}+\sum_{l \in \mathcal{L}}A_{zl}y^{L}_{lt}+\sum_{f \in \mathcal{F}}y^{F}_{fzt}+y^{LShed}_{zt}+\sum_{r \in \mathcal{R}}y^{R}_{rzt} = } & \notag\\
    &&&{\phantom{abc} Y^{D}_{zt}+\sum_{s \in \mathcal{S}}(y^{S+}_{szt}-y^{S-}_{szt})+\sum_{e \in \mathcal{E}}y^{E}_{ezt}+y^{GShed}_{zt}, }  & z \in \mathcal{Z},t \in \mathcal{T}\label{balance}\\
    &\quad&&l^{S}_{sz(t+1)}-l^{S}_{szt}=\eta^{S}_{s}y^{S+}_{szt}-y^{S-}_{szt}, & s \in \mathcal{S}, z \in \mathcal{Z}, t \in \mathcal{T}\label{gmod:storage_balance}\\
    & && \sum_{e \in \mathcal{E}}y^{E}_{ezt}+\sum_{s \in \mathcal{S}^{Hy}} \eta^{EF}v^{SHy-}_{szt} =& \notag\\
    &&&\phantom{abc}\sum_{s \in \mathcal{S}^{Hy}}\eta^{ES}v^{SHy+}_{szt}+\sum_{f \in \mathcal{F}} \eta^{EF}y^{F}_{fzt}, & z \in \mathcal{Z}, t \in \mathcal{T}\label{gmod:hydrogen_balance}\\
    &&&\sum_{s \in S^{T}}\sum_{z \in \mathcal{Z}}\sum_{t \in \mathcal{T}}\sum_{g \in \mathcal{G}}W_{s}\rho^{E}_{g}y^{G}_{gzt} \leq {E}^{\text{CO}_2}, \label{gmod:emission}\\
    &&&y^{P}_{pz}, y^{G}_{gzt}, y^{F}_{fzt}, y^{LShed}_{zt}, y^{R}_{rzt}, y^{S+}_{szt}, y^{S-}_{szt}, y^{E}_{ezt} \in \mathbb{R}^{+}_{0}, \label{gmod:domain}\notag\\
    &&&y^{Gshed}_{zt}, L^{S}_{szt}, v^{SHy-}_{szt}, v^{SHy+}_{szt}, x^{PAcc}_{pz}, x^{LAcc}_{l} \in \mathbb{R}^{+}_{0}.
\end{alignat}
\end{subequations}

The operational cost function $g(x, c)$, which is included in the objective function Equation \eqref{MP_objective} for each operational year $i$, is described by Equation \eqref{SP_objective} that includes total operating costs of all devices, energy load shedding costs and CO$_2$ emissions costs. Vectors $x$ and $c$ contain capacities and costs information, respectively. Constraints \eqref{gmod:tech_cap} and \eqref{gmod:line_cap} ensure devices $(y^{P}_{pzt})$ and transmission lines $(y^{L}_{lt})$ are within their capacities $(x^{PAcc}_{pz}, x^{LAcc}_{l})$. Constraint \eqref{balance} gives the energy balance at each region, where $y^{G}_{gt}, y^{F}_{fzt} \text{ and } y^{R}_{rzt}$ are power generation of generators, fuel cells and renewables respectively. Moreover, we define $y^{E}_{ezt}$ to be the power that goes into electrolysers and $l^{S}_{szt}, y^{S+}_{szt}$ and $y^{S-}_{szt}$ represent the storage level, input and output energy of storage facilities. The energy demand $Y^{D}_{zt}$ can be modelled corresponding to the specific sector, such as offshore platforms. The modelling of offshore platforms is described in details in \ref{full_operational_model}. Constraint \eqref{gmod:storage_balance} states the storage balance of storage facilities. Constraint \eqref{gmod:hydrogen_balance} gives the hydrogen nodal balance of OEHs, where $v^{SHy-}_{szt}$ and $ v^{SHy+}_{szt}$ are the hydrogen output and input of hydrogen storage facilities. Constraint \eqref{gmod:emission} restricts the total emissions. We only consider emissions from the generators, but the model can easily be extend to include other emissions. The complete MILP problem consists of Equations \eqref{MP_objective}-\eqref{mod:operational}.
\section{Results}
\label{results}
We demonstrate the results of a static integrated investment planning and operational problem given by Equations \eqref{MP_objective}-\eqref{mod:operational} for a future point in time. The problem consists of $461,108$ continuous variables, $100$ integer variables and $980,013$ constraints. We implemented the model in Julia $1.6.1$ using JuMP \cite{jump} and solved with Gurobi $9.1.2$ \cite{gurobi}. We run the code on a MacBook Pro with 2.4 GHz 8-core Intel Core i9 processor, with 64 GB of RAM, running on macOS 11.6 Big Sur. The Julia code and data for the case study have been made publicly available \cite{offmod}. We solve the integrated investment and operational model given by Equations \eqref{MP_objective}-\eqref{mod:operational} to conduct sensitivity analysis on CO$_2$ tax, CO$_2$ budget and the capacity of PFS.
\subsection{Case study}
\label{case_study}
The case study is carried out on the North Sea part of the NCS, considering $66$ fields. The problem consists of $77$ regions, divided into $66$ fields, $5$ OHEs and $5$ onshore buses. By using the clustering approach described in Section \ref{modelling_strategies_and_assumption} we find that the system can be represented using $5$ clusters and henceforth go from $77$ regions to $15$ regions. We assume that the power demand of platforms initially is entirely supplied by gas turbines, as only a limited number of platforms receives PFS. We construct four representative months with hourly resolution and scale them up to represent a whole year.

% The oil and gas industry is capital intensive, and it involves many stakeholders. Therefore, 
The operational data of oil and gas industry is sensitive and usually not disclosed to the public. Aggregated data such as monthly or yearly production of petroleum on the NCS can be obtained from \cite{NPD}. One can also find monthly production and injection data for each field from \cite{field_monthly_data} and \cite{diskos}. Neither of these can be directly used as inputs for this study due to the time resolution difference. 
% It is difficult to obtain all the data directly from operators. 
Therefore, reasonable data generation is necessary. Raw data is collected from (a) Norne ($1998-2006$) and Volve ($2008-2016$) fields with hourly production and injection data from \cite{lowemission_centre} and (b) monthly production and injection data of each field from \cite{field_monthly_data}. We come up with a data generation method that considers the lifetimes of offshore fields
\cite{offmod}.

We define a base case (Base) with offshore renewables, electric boiler, battery and PFS as investment options. We then use this case as a benchmark to check against the case with OEHs. The full model given by Equations \eqref{MP_objective}-\eqref{mod:operational} takes approximately 2 hours.
% =============================================================================
% results
% =============================================================================
\subsection{Energy system analysis}
\label{enery_systme_analysis}
\begin{figure}[!thb]
    \centering
    \includegraphics{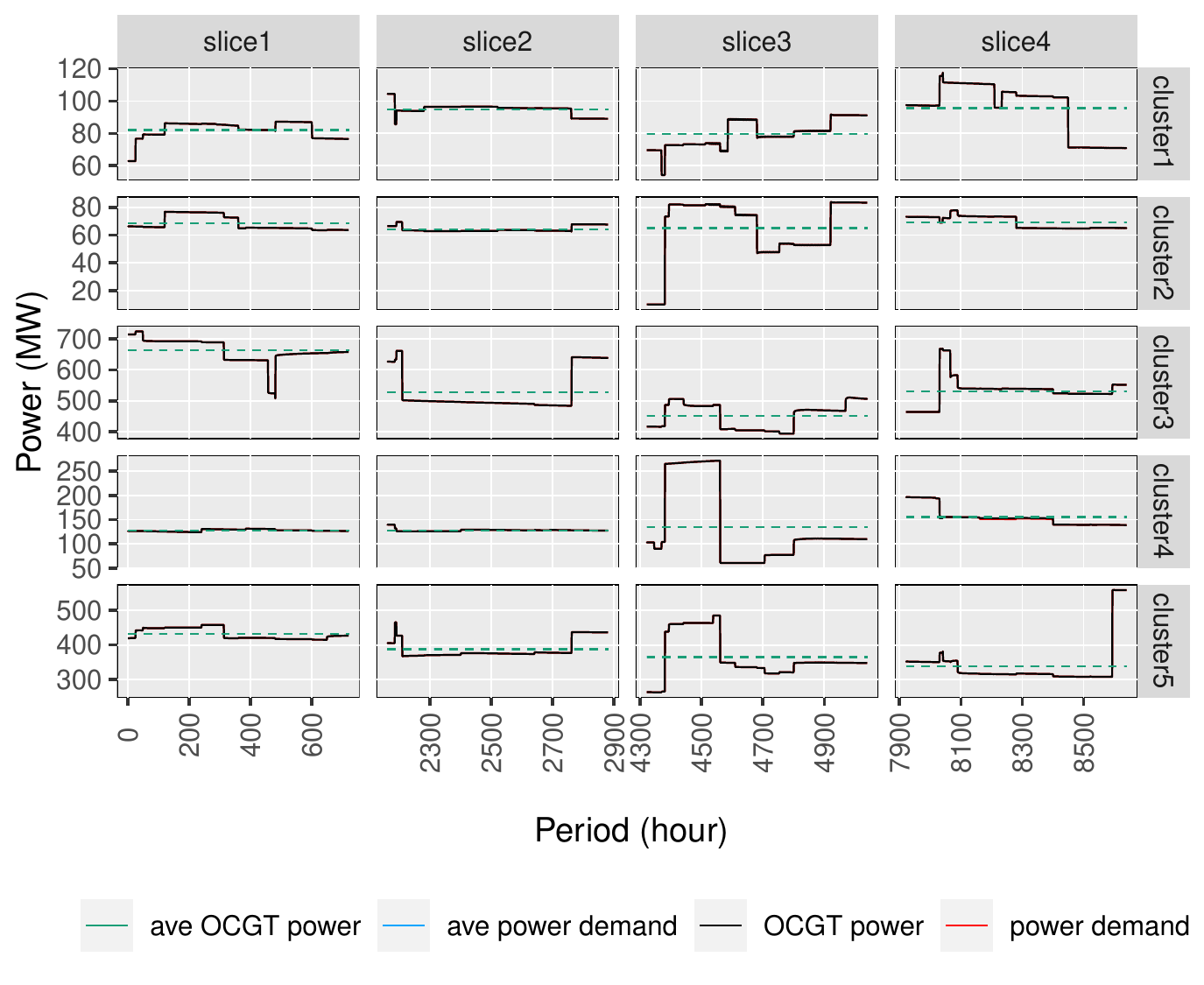}
    \caption{Power consumption and supply (Only two lines are observable since power supply and demand match exactly. OCGT power equals power demand at all times).}
    \label{electricity_load_supply}
\end{figure}
\begin{figure}[!thb]
    \centering
    \includegraphics{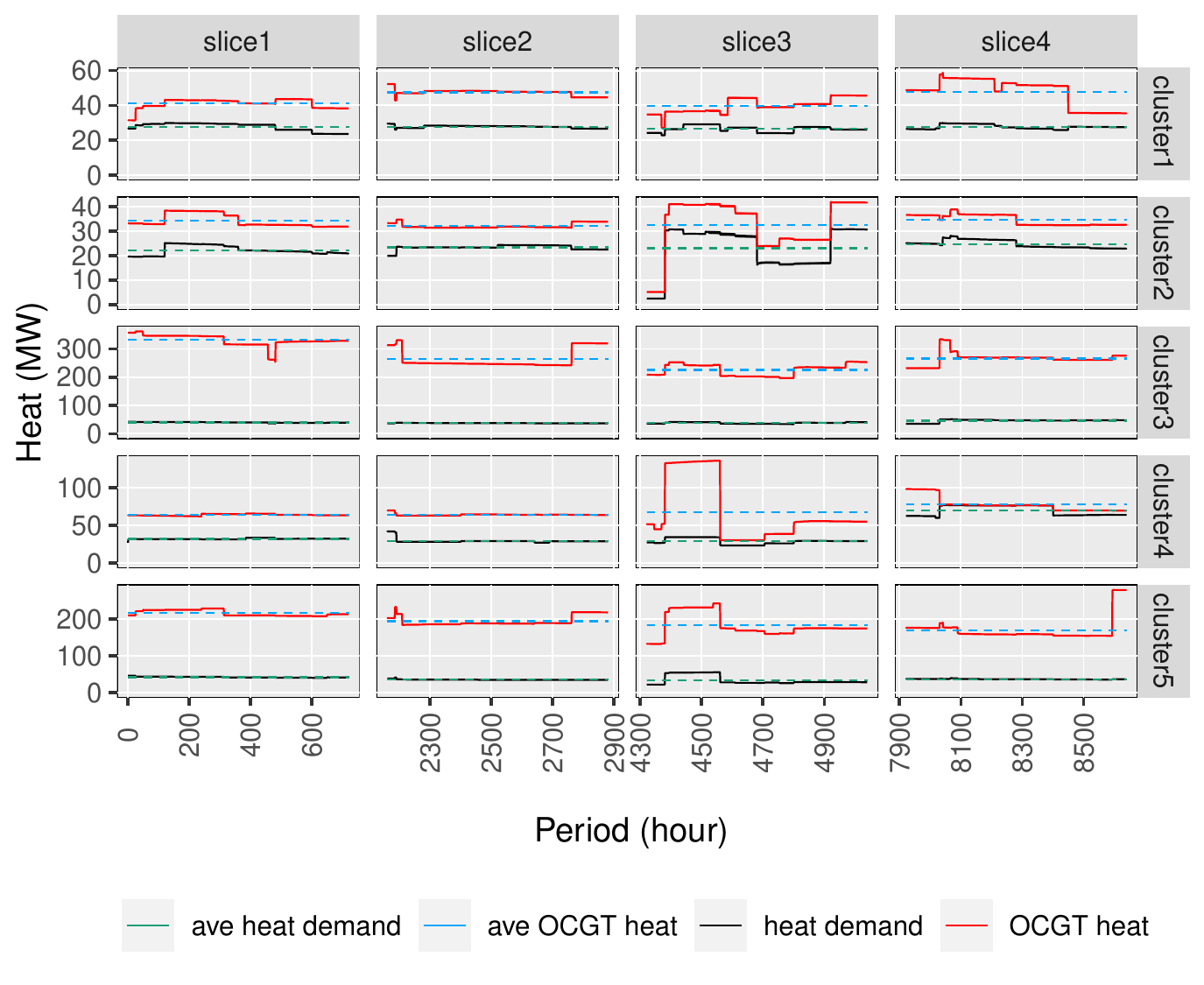}
    \caption{Heat consumption and supply.}
    \label{heat_load_supply}
\end{figure}
In this section, we present results on energy consumption and CO$_2$ emission of the initial system. By post-processing, we verify that the energy consumption of platforms is of the same order of magnitude as reported numbers. The resulted CO$_2$ emission is $5.51$ Mt/yr. In comparison, the reported total emission of the relevant fields was $6.89$ Mt in $2019$ \cite{diskos}. We expect emissions from the model to be lower than $6.89$ Mt since not all emission sources are considered. Based on \cite{nguyen}, one could assume that the major processes considered in this study make up about $80\%$ of the total load. Therefore, $5.51$ Mt yearly emission is within the correct range, implying that the energy load modelling is relatively accurate.

\begin{figure}[!htb]
    \centering
    \includegraphics{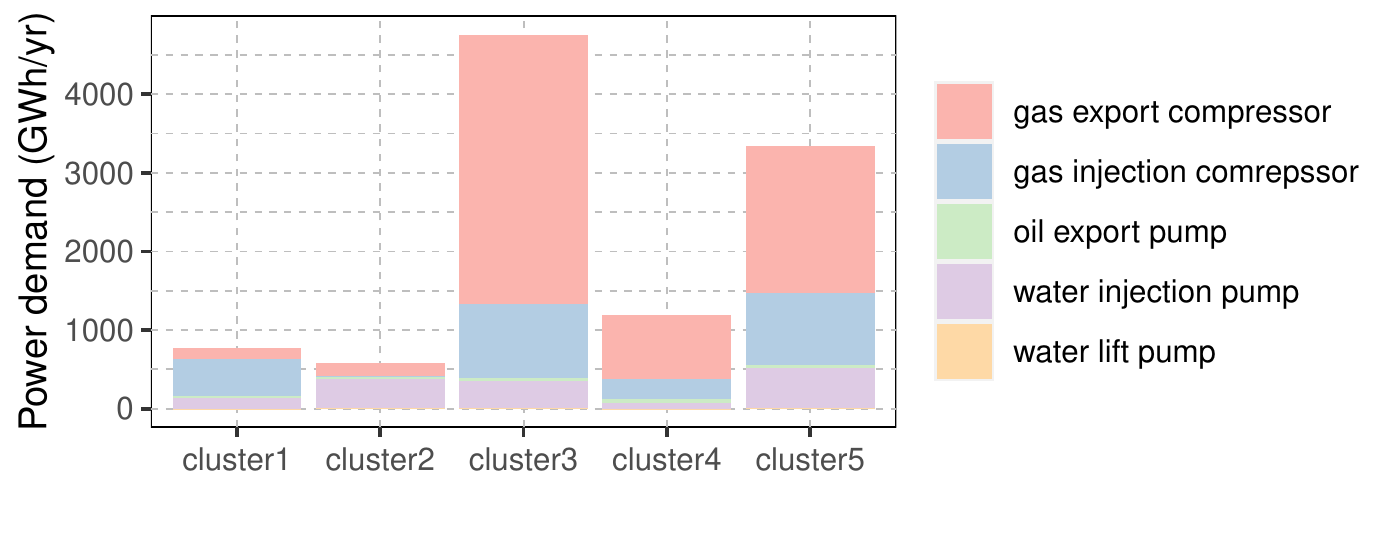}
    \caption{Power demand in a year.}
    \label{rel_power_demand}
\end{figure}
% \begin{figure}[!htb]
%     \centering
%     \includegraphics{rel_emission_pie.pdf}
%     \caption{Emission distribution by cluster.}
%     \label{rel_emission}
% \end{figure}
\begin{table}[!]
\begin{tabular}{cccccc}
\hline
& cluster1 & cluster2 & cluster3 & cluster4 & cluster5 \\ \hline
Emission distribution & 7.2\%    & 5.5\%    & 44.7\%  & 11.2\%   & 31.3\%   \\ \hline
\end{tabular}
\caption{Emission distribution by cluster.}
\label{table:rel_emission}
\end{table}
From Figure \ref{electricity_load_supply}, we can see that the power output of the Open Cycle Gas Turbine (OCGT) matches power demand at every operational period. Heat recovery of OCGTs is assumed the only heat source. Figure \ref{heat_load_supply} shows that heat recovery of OCGTs provides more than enough heat due to high electricity generation.
% and there is no need for OCGT to generate extra electricity for supplying heat rather than power. 
We can also see that energy consumption can vary significantly. A breakdown of electricity load is shown in Figure \ref{rel_power_demand}, gas export compressors dominate the power consumption in clusters $3-5$. Water injection is the largest power consumer in cluster $2$ since there are some mature fields (e.g., Ekofisk) whose reservoir pressures are mainly maintained by water injection. OCGT is the only energy and emission source in the initial setup. Therefore, emission breakdown includes the emissions from the total energy consumption of each region. Cluster $1$ has the second smallest share of the total energy consumption, with a considerable amount of power consumed by gas injection. The fields in cluster $1$, such as Grane, have the third-highest gas injection level among the $66$ fields. From Table \ref{table:rel_emission}, we find that emission mainly comes from the northern part of the North Sea.
\subsection{Sensitivity analysis of \texorpdfstring{CO$_2$}{} tax}
This section presents the results of sensitivity analysis of CO$_2$ tax. We introduce CO$_2$ tax and still keep the carbon budget inactive. We increase the carbon tax from $55$ to $500$ €/tonne with a step size of $5$ €/tonne. PFS capacity limits are estimated from \cite{power_from_land,esurplus}.  Note that the cost of PFS may be underestimated since we only consider the costs of subsea cables, onshore and offshore converter stations and electricity bills. In reality, PFS projects may also face investment in onshore transmission lines or onshore power system capacity expansion. 
% AACs may underestimate the actual costs of PFS, especially when PFS is used for compensating for offshore wind power. However, the underestimation can be mitigated by adjusting AACs via post-processing. More specifically, the capacities of PFS from the model output should be aligned with the PFS capacities reported that are correspondent to the AACs.  We use the maximum PFS as the rated capacity. 
We analyse the results from two metrics, CO$_2$ emission and energy loss. Energy losses are from conversions, transmission, and generation shed. The calculation is presented in \ref{sec:energy_loss}.
\begin{figure}[t]
    \centering
    \includegraphics{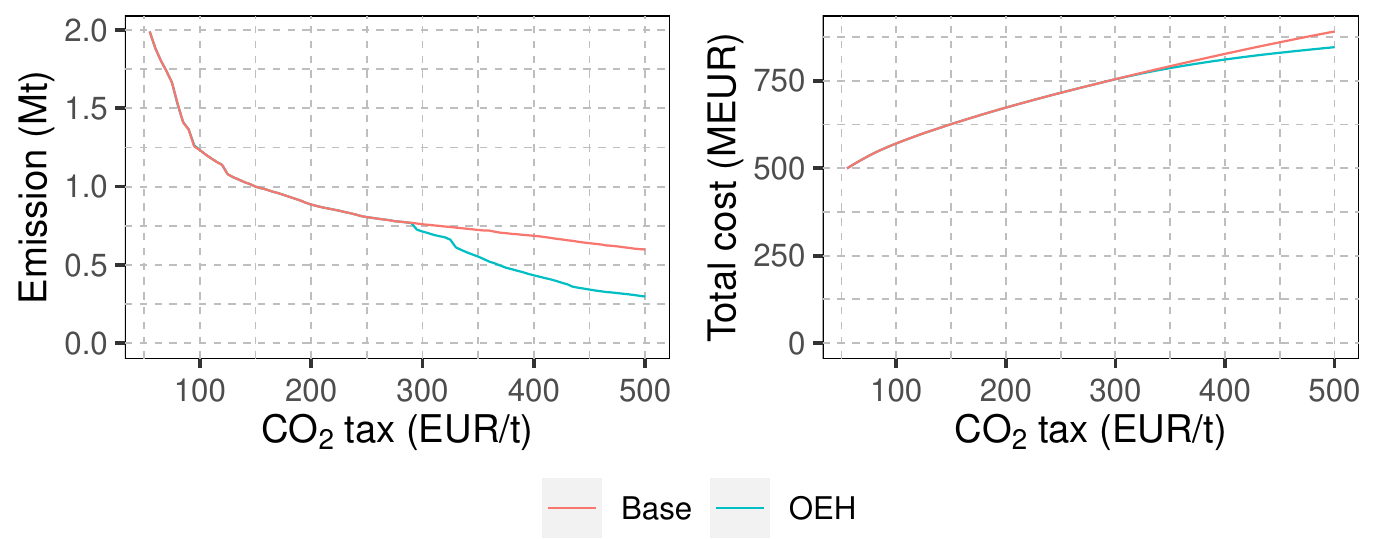}
    \caption{Emission and cost comparison (CO$_2$ tax sensitivity analysis).}
    \label{co2tax_emission_cost}
\end{figure}
\begin{figure}[!]
    \centering
    \includegraphics{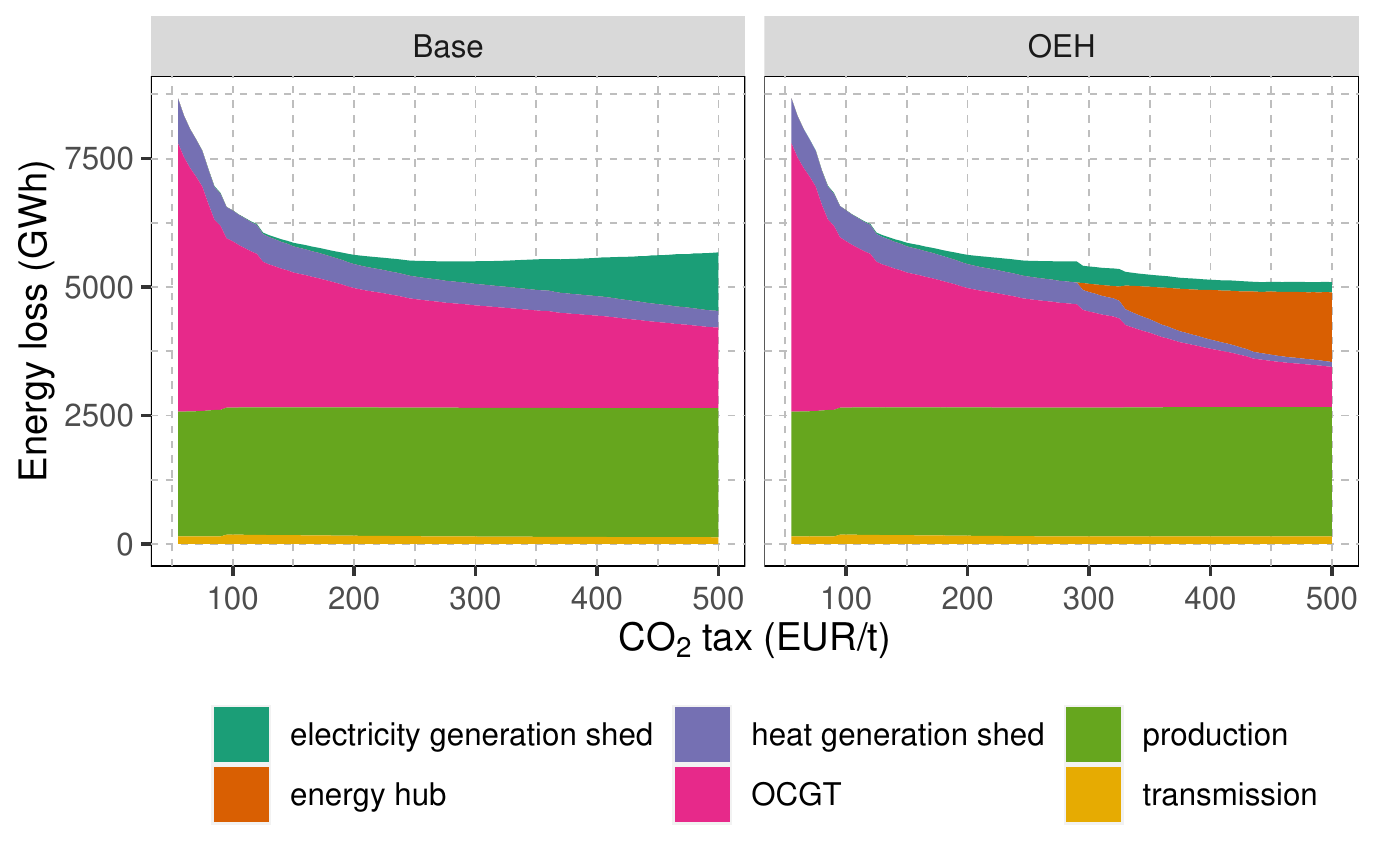}
    \caption{Energy loss (CO$_2$ tax sensitivity analysis).}
    \label{co2tax_energy_loss}
% \end{figure}
% \begin{figure}[!htb]
    \centering
    \includegraphics{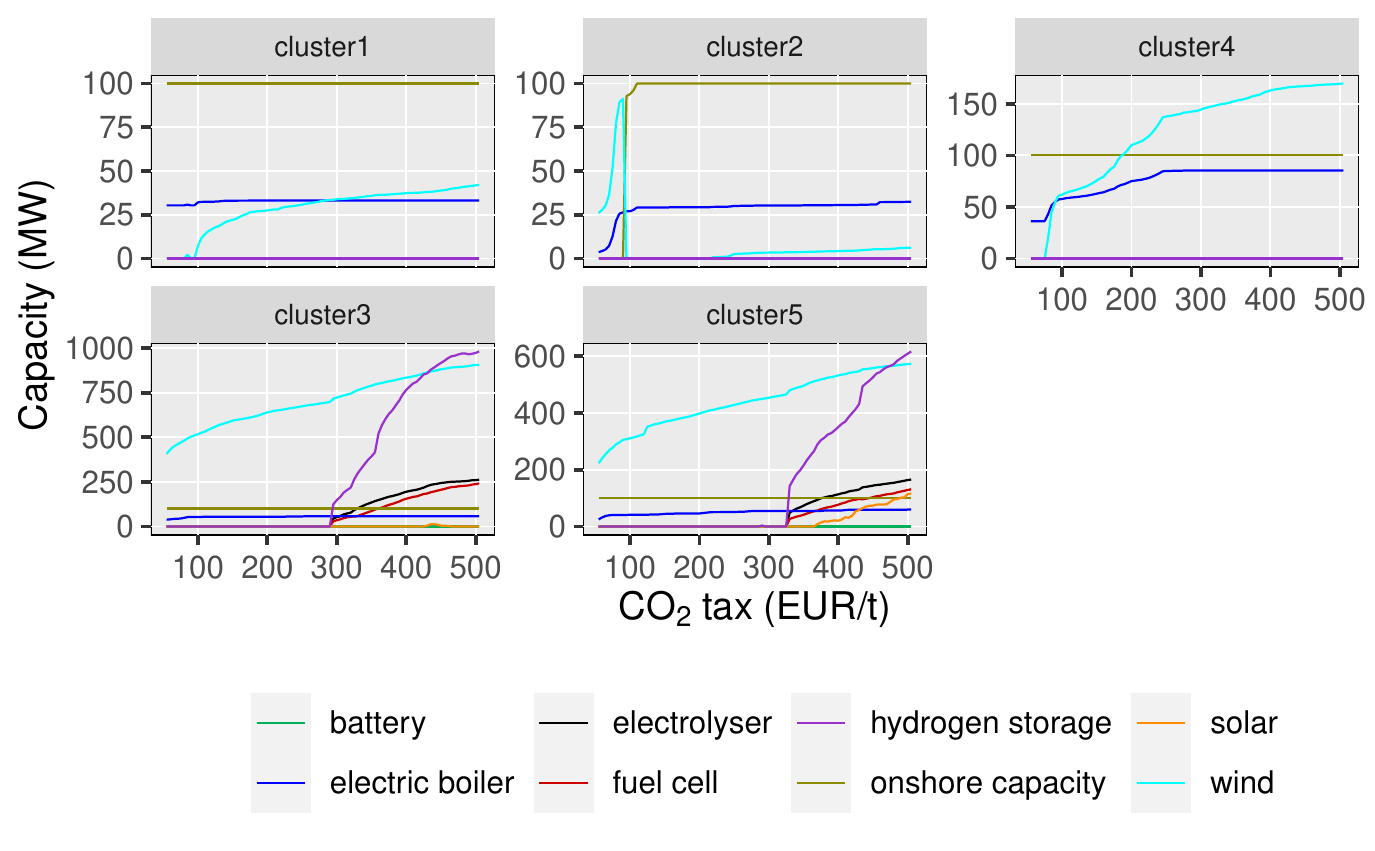}
    \caption{Capacities of technologies in each cluster (CO$_2$ tax sensitivity analysis),  hydrogen storage is measured in tonne.}
    \label{co2tax_tech_cap}
\end{figure}
% \begin{table}[t]
%     \caption{Combination of technologies under different CO$_2$ tax levels (€/tonne, numbers indicate the rage of CO$_2$ tax).}
%     \centering
%     \resizebox{\textwidth}{!}{
% \begin{tabular}{llllllllll}
% \hline
% cluster & OCGT & \begin{tabular}[c]{@{}l@{}}OCGT\\ electric boiler\end{tabular} & \begin{tabular}[c]{@{}l@{}}OCGT\\ wind\end{tabular} & \begin{tabular}[c]{@{}l@{}}OCGT\\ wind\\ electric boiler\end{tabular} & \begin{tabular}[c]{@{}l@{}}OCGT\\ wind\\ electric boiler\\ energy hub\end{tabular} & \begin{tabular}[c]{@{}l@{}}OCGT\\ wind\\ electric boiler\\ PFS\end{tabular} & \begin{tabular}[c]{@{}l@{}}OCGT\\ wind\\ electric boiler\\ energy hub\\ PFS\end{tabular} & \begin{tabular}[c]{@{}l@{}}OCGT\\ wind\\ solar\\ electric boiler\\ energy hub\end{tabular} & \begin{tabular}[c]{@{}l@{}}OCGT\\ wind\\ solar\\ electric boiler\\ energy hub\\ PFS\end{tabular} \\ \hline
% cluster1 &0 -- 32      &-              &34 -- 38     &40 -- 288    &290 -- 352   &-   &354 -- 500  &- &- \\
% cluster2 &0 -- 28      &-              &30 -- 34     &36 -- 310    &312 -- 452   &-   &454 -- 500  &- &-\\
% cluster3 &0 -- 36      &-              &38 -- 40     &42 -- 152    &-   &154 -- 254  &256 -- 386, 390-394  &- &388, 396 -- 500\\
% cluster4 &-             &0 -- 34       &-              &36 -- 252    &-   &254 -- 374    &376 -- 500    &- &- \\
% cluster5 &0 --34       &-              &36 -- 40     &42 -- 248    &250 -- 328   &-   &- &330 -- 452 &454 -- 500\\ \hline
% \end{tabular}
% }
% \label{co2_tax_tech_table}
% \end{table}

From Figure \ref{co2tax_emission_cost}, we can see that CO$_2$ tax as a single instrument may not be enough to motivate a zero emission system. We also find that near zero emission can be achieved with a very high CO$_2$ tax. Therefore, a hard carbon cap may be necessary for stimulating a zero emission system. When CO$_2$ tax is $55$ €/tonne, the system reduces about 64\% of the initial emissions. Approximately $84\%$ of the emissions can be cut if CO$_2$ tax is increased to $200$ €/tonne as planned. As OCGTs are replaced by renewable energy, energy loss is reduced as well. OEHs can potentially reduce up to around $50\%$ more CO$_2$ emission and $5\%$ total cost than the case with only offshore wind and PFS (Base) at certain CO$_2$ tax levels. From Figure \ref{co2tax_energy_loss}, we find that energy loss during production counts for $11\%$ of the energy loss. OCGTs lose $18$ GWh of energy during an operational year.
As production from wind turbines replaces gas turbines, energy loss from OCGT is reduced. However, due to the lack of energy storage, electricity generation shedding increases because wind power is shed. We find that OEHs can effectively reduce electricity generation shedding although, it loses energy during conversion. Overall, energy loss is up to $10\%$ lower in the case of OEHs compared with Base at certain tax levels.

From Figure \ref{co2tax_tech_cap},
we find that different clusters show different levels of sensitivity to CO$_2$ tax. Offshore wind is the first renewable energy solution that is added to the system. Electric boilers are soon needed as offshore wind replaces gas turbines partially. OEHs are installed when CO$_2$ tax is above $290$ €/tonne. Offshore solar is only added in clusters $3$ and $5$ under very high CO$_2$ tax levels. OCGTs still operate even CO$_2$ tax increases to $500$ €/tonne.
% Figure \ref{co2tax_tech_cap} shows the capacities in each technologies in detail. 
We can see that in a static planning problem, if CO$_2$ tax is the only instrument and increases to $200$ €/tonne as the government's plan in 2030, OEHs may not be necessary. However, CO$_2$ tax combined with the EU emissions trading system may likely increase the total CO$_2$ price to around $250-300$ €/tonne, which is about the breakeven price of OEHs. In addition, the potential benefits of the OEHs may realise once they provide services to more sectors, such as exporting hydrogen for industries or transportation.
% We also notice that using MILP model gives more realistic investment decisions. For example, linear model would make investment in HVDC cables with several megawatts capacities which may be unrealistic and doesn't exist in a mixted-intger linear model.
%%%%%%%%%%%%%%%%%%%%%%
\subsection{Sensitivity analysis of \texorpdfstring{CO$_2$}{} budget}
\begin{figure}[t]
    \centering
    \includegraphics{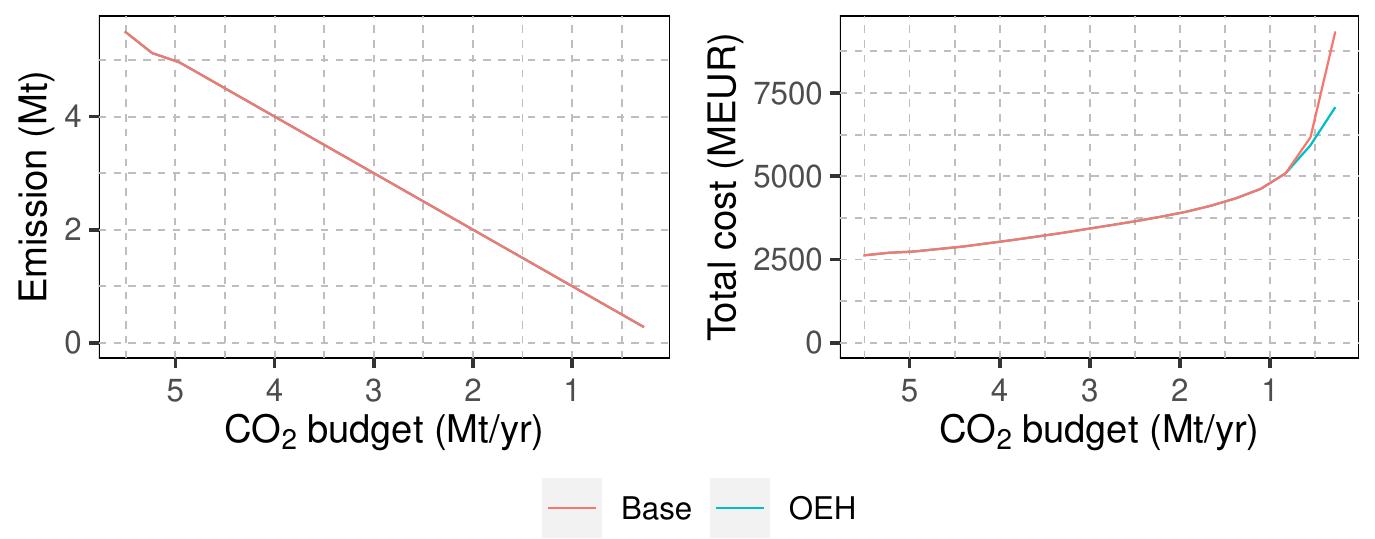}
    \caption{Emission and cost comparison (CO$_2$ budget sensitivity analysis).}
    \label{co2budget_emission_cost}
\end{figure}
\begin{figure}[!]
    \centering
    \includegraphics{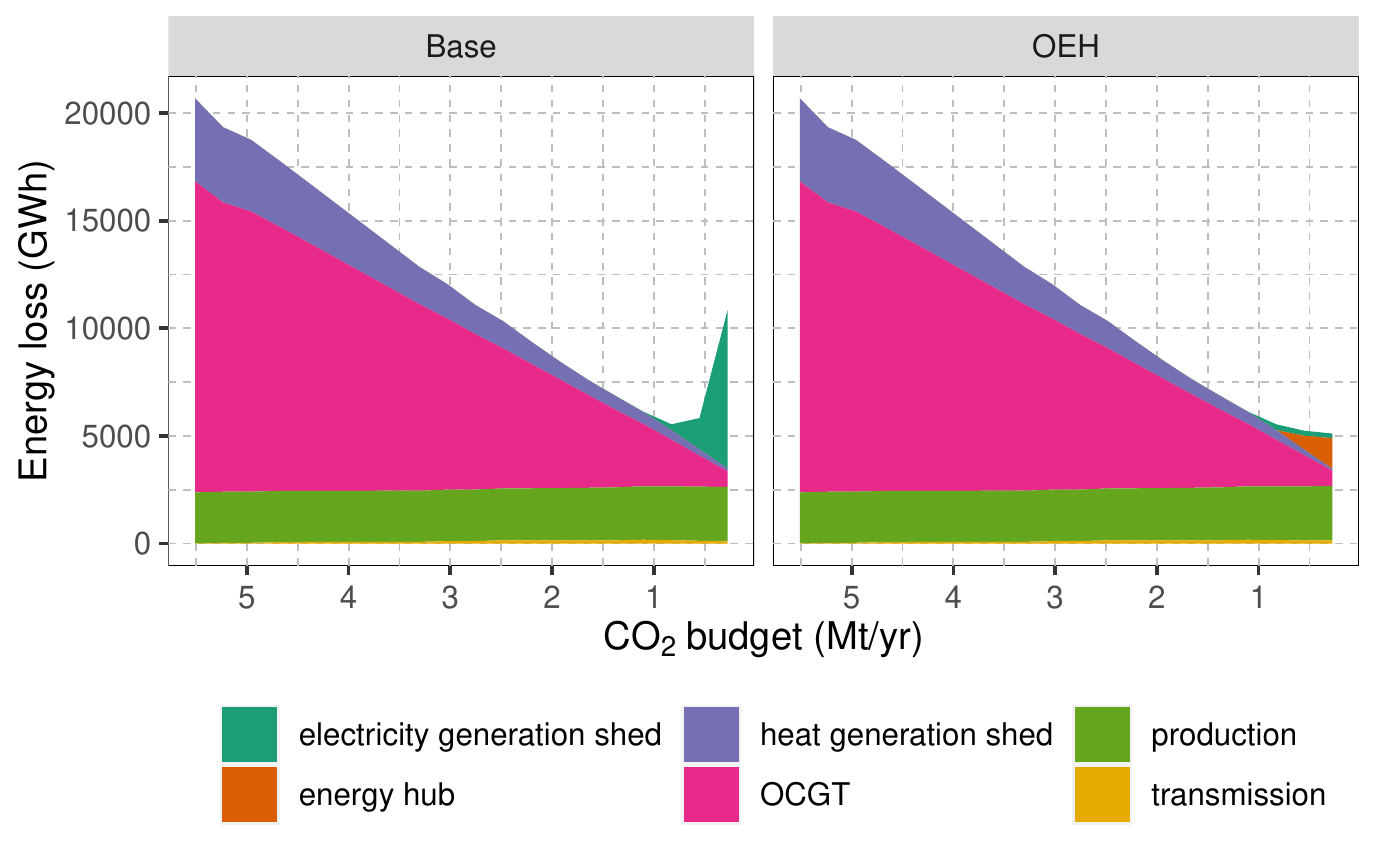}
    \caption{Energy loss (CO$_2$ budget sensitivity analysis)}
    \label{co2budget_energy_loss}
% \end{figure}
% \begin{figure}[!]
%     \centering
    \includegraphics{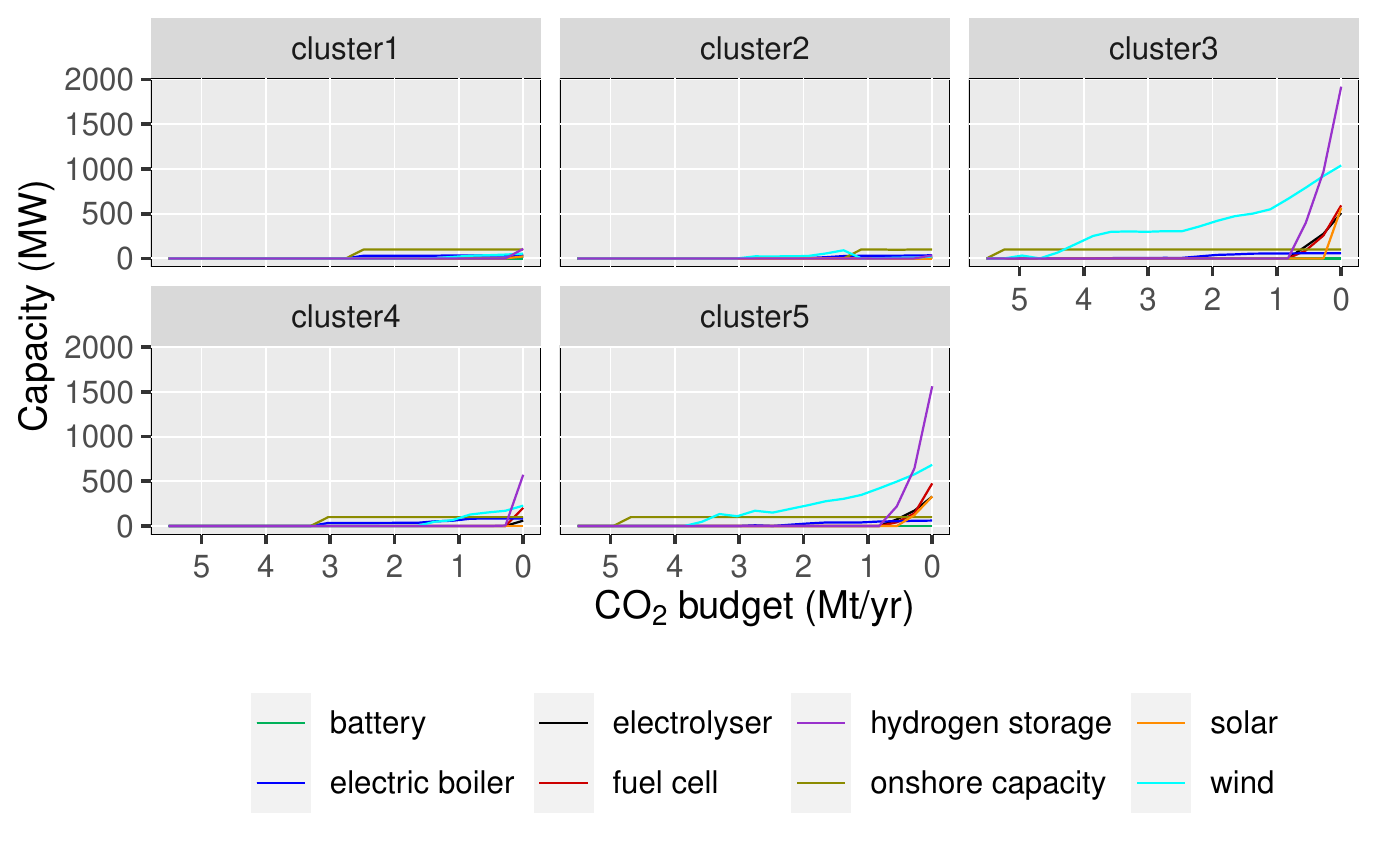}
    \caption{Capacities of technologies in each cluster (CO$_2$ budget sensitivity analysis),  hydrogen storage is measured in tonne.}
    \label{co2budget_tech_cap}
\end{figure}
For the CO$_2$ budget, we use initial emission as the starting point and reduce it by $5\%$ until it hits $0$. From Figure \ref{co2budget_emission_cost}, we find that carbon cap is binding most of the time and we rarely see that emissions reduce more than the carbon cap. Thus, there is no difference in actual emission in Base and the system with OEHs. However, the cost is 25\% lower in a zero emission system with OEHs compared with Base.
% CO$_2$ tax is fixed to the current level, $55$ €/tonne.The capacity limits and AACs of PFS remain the same. 

We find that in a zero emission system without OEHs, energy loss is around 477 TWh due to 81 GW of wind power capacity without storage. This may not be likely from happening since some forms of storage would be added to compensate for offshore wind in reality. From Figure \ref{co2budget_energy_loss}, we can see a large amount of energy loss when reaching near zero emission system in Base. The energy loss in Base is 10,832 GWh in a near zero emission system that is about twice higher than that with OEHs. A large amount of wind power is installed to meet power demand at any time. Therefore, the same capacity of wind that can cope with peak demand hours will also generate surplus power during normal hours. This leads to increased energy losses as more wind replaces OCGT without proper energy storage. In the case of OEHs, wind power can be stored when excess power is generated. It is also worth noticing that in the energy system without an OEH, energy storage is the battery on the platforms, which can be infeasible due to space and weight limitations. We observe that investments in batteries are only needed when approaching zero emission in Base. No battery is needed in a system with OEHs. In addition, the energy loss of OEHs is 27\% of the total loss, and the loss during production is about 50\% of the total. 

From Figure \ref{co2budget_tech_cap}, we find that cluster 3 receives PFS after a 10\% reduction of the carbon cap. Cluster 3 has the highest emission level but the shortest distance from shore. Therefore, taking PFS and partially electrify the fields in cluster 3 can help the system reduce 10\% emission in a cost effective way. The system does not cut emissions proportionally in each cluster but cut emissions from clusters with the highest emission, such as cluster 3 and cluster 5. Therefore, it may be necessary to consider the whole NCS when conducting system planning rather than consider each cluster separately and reach sub-optimality. Cluster 2 is the most remote, more than 300 km from the shore; PFS is less economical than offshore wind. Therefore, offshore wind is added to cluster 2 when the carbon cap drops to 2.75 Mt/yr. 
% we find that offshore wind and PFS can reduce emissions and help the system align with the CO$_2$ budget until the carbon cap reduces to around $0.55$ Mt/yr. After that, either an OEH or PFS is needed for compensating for wind volatility.
When the CO$_2$ budget reduces to below $0.83$ Mt/yr, CO$_2$ emissions are nearly zero in clusters $1$ and $2$. However, the carbon cap needs to reduce to zero to shut down OCGTs completely in all clusters. Nearly $4,200$ tonnes of hydrogen storage capacity is needed in a zero emission NCS energy system, and nearly half is installed in cluster 3. 
%
% Figure \ref{co2budget_co2_marginal_price} shows a comparison of CO$_2$ marginal price in two cases. We find that the CO$_2$ marginal price can potentially be ten times higher in the case without OEH. This implies that OEH can potentially provide more cost-efficient decarbonisation solutions.
% The mixtures of technologies when achieving zero emission are summarised in Table .
\subsection{Sensitivity analysis of the capacity of PFS}
We now present the results of sensitivity analysis of the capacity of PFS. The capacity of PFS affects the investments in offshore technologies. An onshore system has a limited capacity to transmit power offshore. Furthermore, onshore system expansion can vary the limitation. However, we do not consider onshore system expansion in this paper. Therefore, we conduct sensitivity analysis to reveal the relationship between onshore power system capacity and offshore decarbonisation technologies.
\subsubsection{Scenario 1 (S1)}
The first scenario is to fix the CO$_2$ tax to $300$ €/tonne and increase the PFS capacity of each onshore node from 0 MW to 1000 MW with a 10 MW step. We find that the investment decisions remain the same if the PFS capacity is higher than 710 MW. Therefore, we only present the results from 0 MW to 710 MW.
% From previous cases, we have seen that wind power alone can reduce emissions to less than 2 Mt/yr. Therefore, we aim to find the relationship between PFS and offshore wind in this scenario.
\begin{figure}[t]
    \centering
    \includegraphics{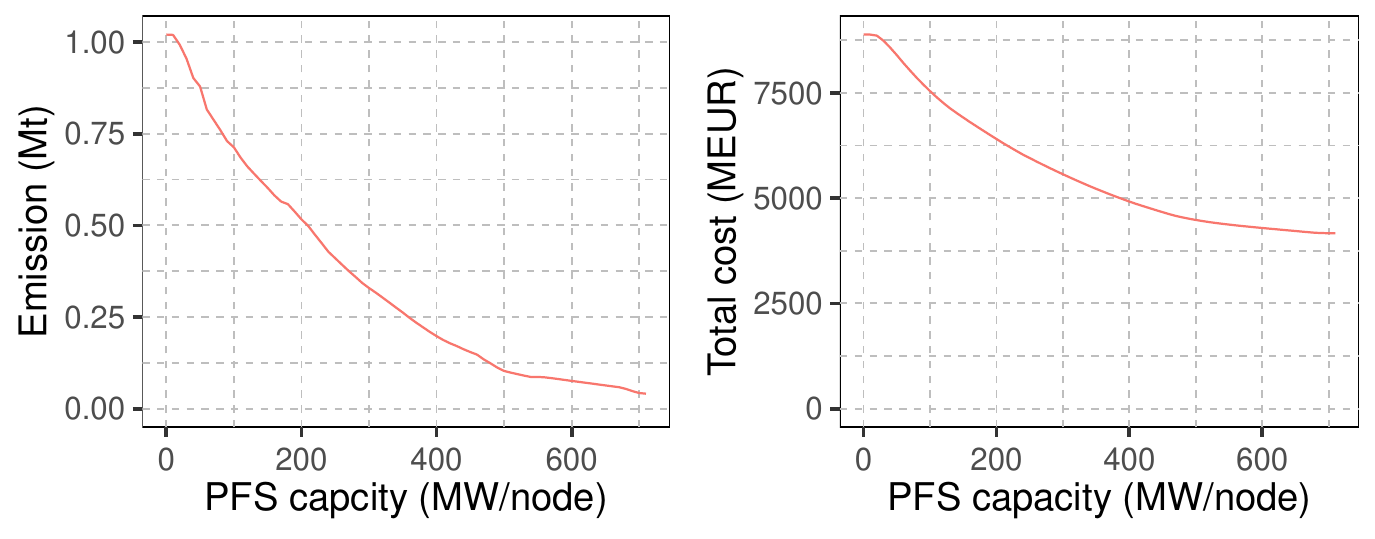}
    \caption{Emission and cost (PFS capacity sensitivity analysis, S1).}
    \label{PFS_emission_cost_s1}
\end{figure}
\begin{figure}[!]
    \centering
    \includegraphics{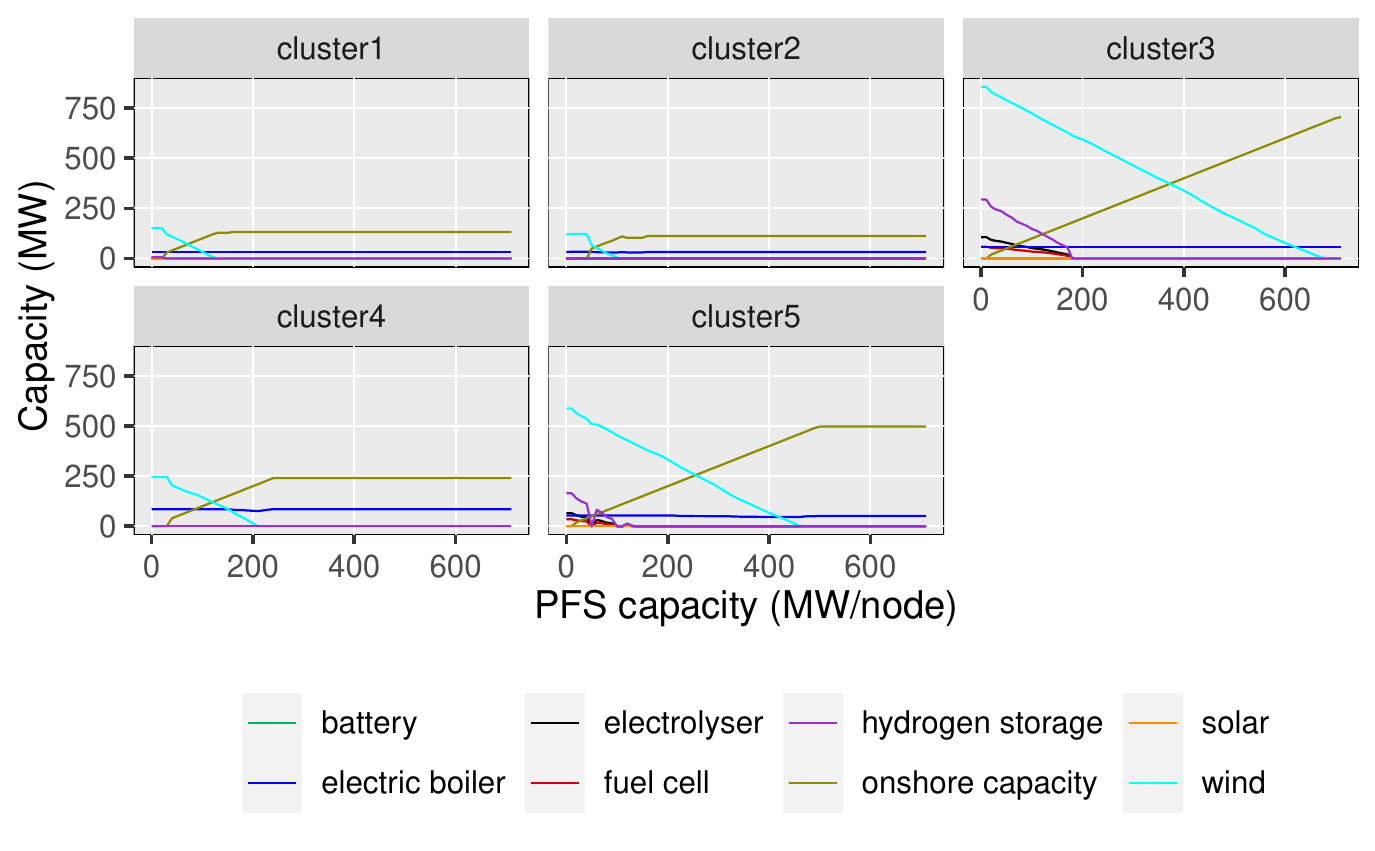}
    \caption{Capacities of technologies in each cluster (PFS capacity sensitivity analysis, S1), hydrogen storage is measured in tonne.}
    \label{PFS_tech_s1}
% \end{figure}
% \begin{figure}[!htb]
%     \centering
    \includegraphics{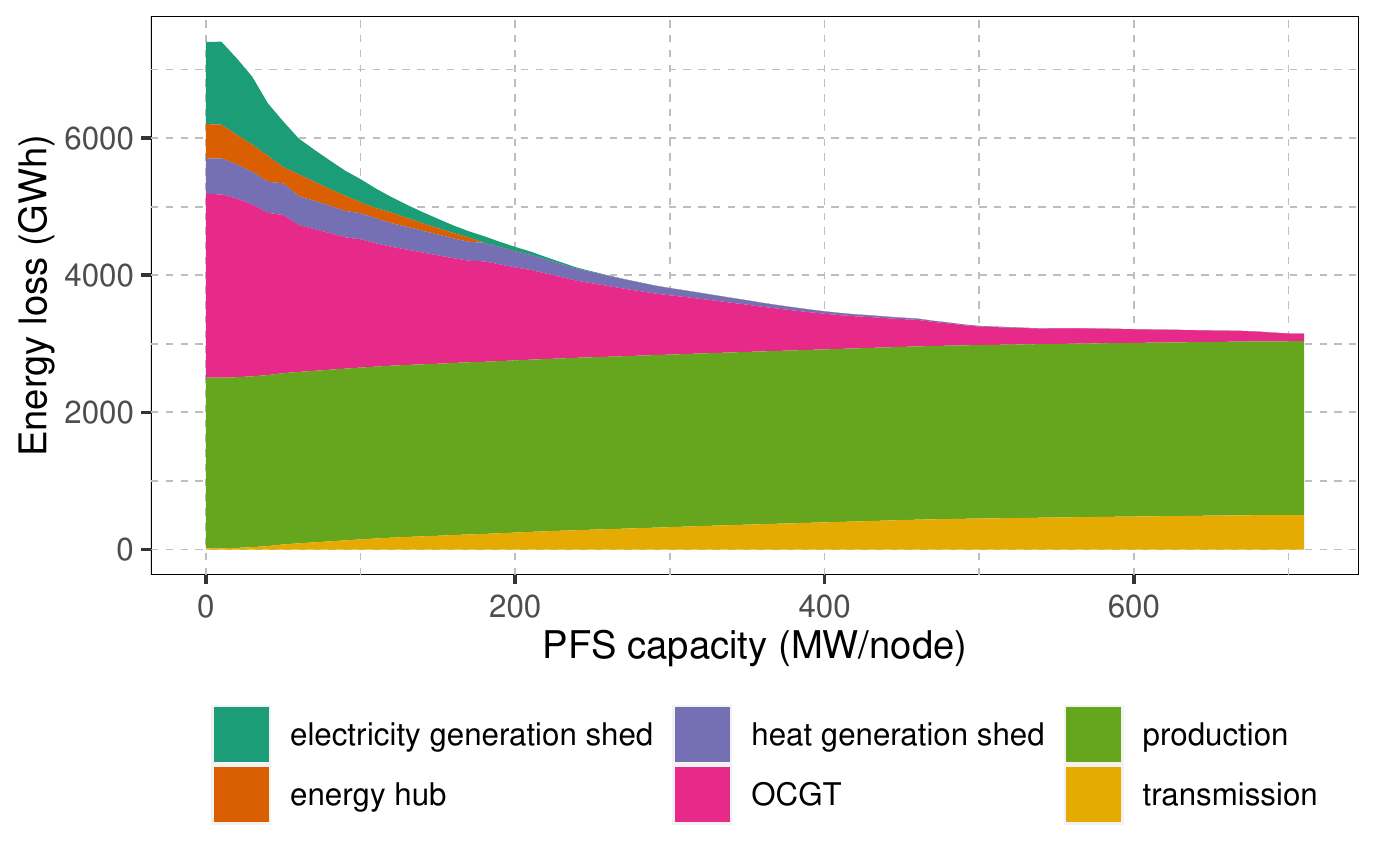}
    \caption{Energy loss (PFS capacity sensitivity analysis, S1).}
    \label{PFS_energy_loss_s1}
\end{figure}
% \begin{table}[t]
%     \label{abatement_krt_tech_table_s1}
%     \caption{Combination of technologies under different ACC of PFS (S1).}
%     \centering
%     % \resizebox{\textwidth}{!}{
%     \begin{adjustbox}{width=\columnwidth,center}
% \begin{tabular}{llllll}
% \hline
% cluster  & \begin{tabular}[c]{@{}l@{}}wind\\ electric boiler\end{tabular} & \begin{tabular}[c]{@{}l@{}}wind\\ PFS\\ electric boiler\end{tabular} & \begin{tabular}[c]{@{}l@{}}wind\\ PFS\end{tabular} & PFS & \begin{tabular}[c]{@{}l@{}}PFS\\ electric boiler\end{tabular} \\ \hline
% cluster1 &104 -- 800 &84 -- 100 &80 &64 -- 76 &4 -- 60 \\
% cluster2 &104 -- 800 &80 -- 100 &76 &64 -- 72 &4 -- 60 \\
% cluster3 &104 -- 800 &88 -- 100 &84 &64 -- 80 &4 -- 60 \\
% cluster4 &104 -- 800 &84 -- 100 &-- &-- &4 -- 80 \\
% cluster5 &104 -- 800 &88 -- 100 &84 &64 -- 80 &4 -- 60 \\ \hline
% \end{tabular}
% % }
% \end{adjustbox}
% \end{table}
From Figure \ref{PFS_emission_cost_s1}, we can see that by having 710 MW capacity in each onshore node, the system can achieve 0.05 Mt/yr emission and reduce nearly 50\% of the total cost. However, increasing the capacity further does not cut emissions or costs further. Figure \ref{PFS_energy_loss_s1} shows that energy loss during transmission makes up 16\% of the total energy loss as we increase the onshore capacity. Electricity generation shed decreases as onshore capacity increases because PFS gradually replaces offshore wind, and less energy is lost from wind turbines. From Figure \ref{PFS_tech_s1}, we find that for onshore nodes that connect to cluster 1 and cluster 2, the needed onshore capacities are about 130 MW and 110 MW, respectively. There are also upper limits on the installed capacity of PFS in the other clusters. We also notice that OEHs are still needed in clusters 1, 3, 5 as we increase the onshore capacity. However, eventually, OEHs are not needed since PFS can provide more stable power and OEH with storage becomes less important.
% The role of the technologies in the energy mix does not change until the cost of OEHs drop \%. Therefore, we only present the segment $4$ €/tonne -- $120$ €/tonne. There are two significant drops showed in Figure \ref{abatement_krt_emission_s1} when AACs are around $100$ €/tonne and $60$ €/tonne separately. The first drop occurs because  $100$ €/tonne is the breakeven point for PFS projects. In addition to wind power, onshore power is also added to the system. Then when AACs drop to around $60$ €/tonne, PFS becomes more cost-effective than offshore wind and replaces the wind power. When the AAC reduces to around $52$ €/tonne, the system achieves near zero emission. Figure \ref{aac_energy_loss_s1} shows that transmission loss counts for $17\%$ of the total energy loss when approaching zero emission. Figure \ref{aac_tech_cap_s1} shows that 1000 MW offshore wind is needed to reduce half the initial emission. When AACs reduce to around $60$ €/tonne, PFS is more cost-efficient than local OCGTs.
\subsubsection{Scenario 2 (S2)}
\begin{figure}[!]
    \centering
    \includegraphics{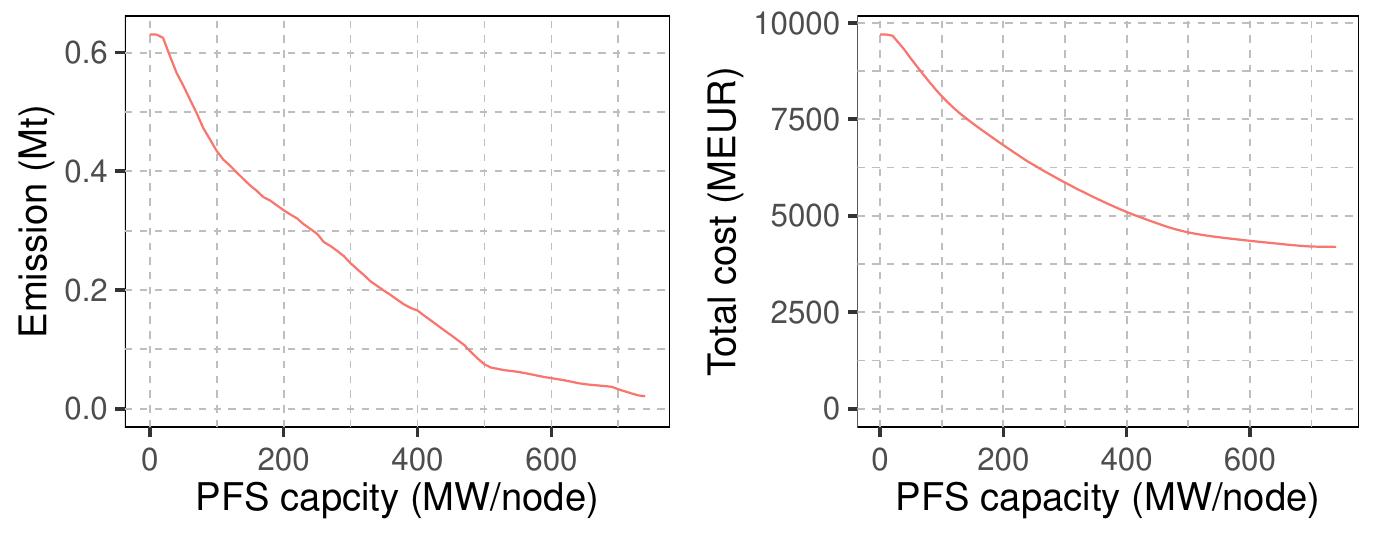}
    \caption{Emission and cost (PFS capacity sensitivity analysis, S2).}
    \label{PFS_emission_cost_s2}
\end{figure}
\begin{figure}[!]
    \centering
    \includegraphics{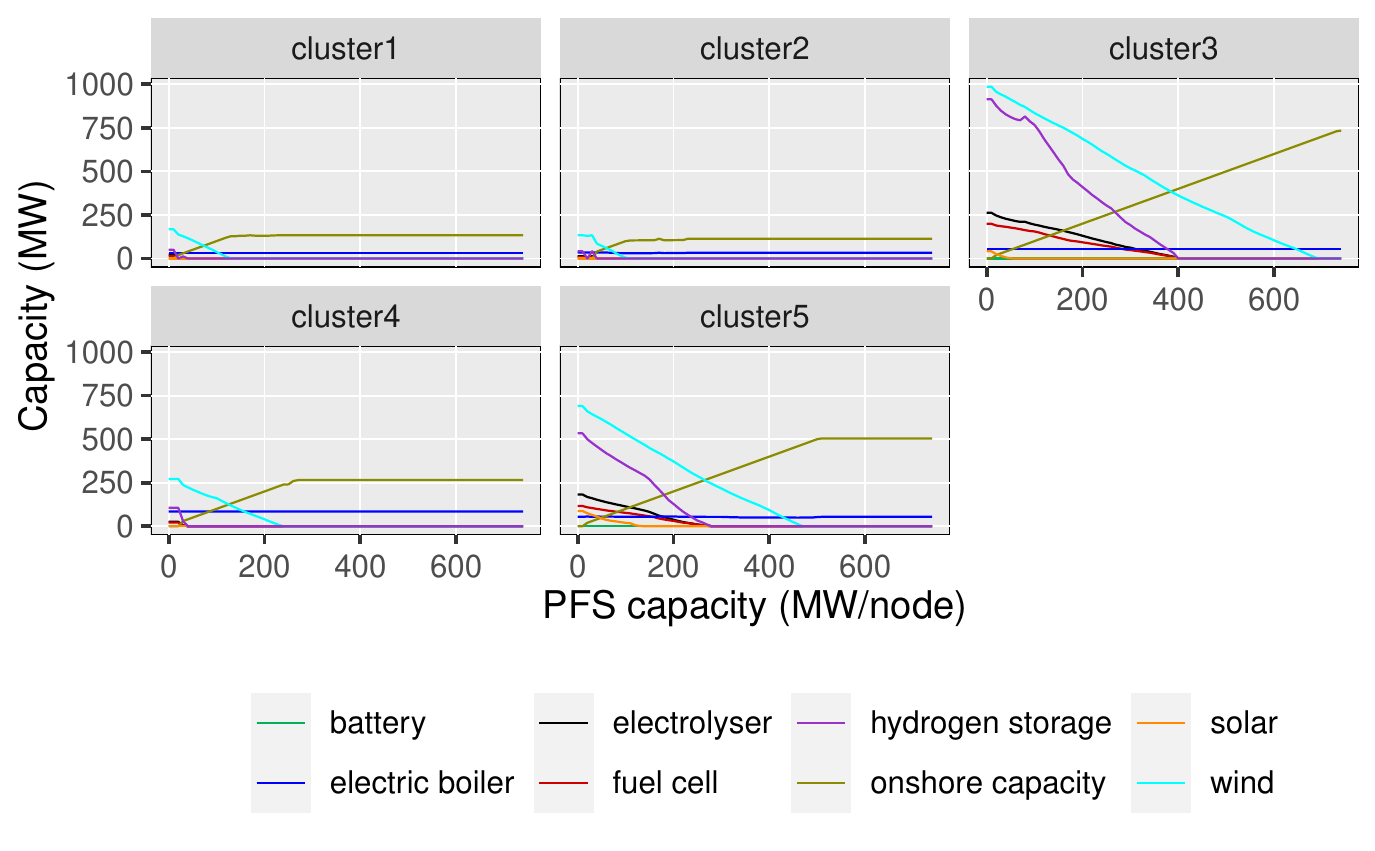}
    \caption{Capacities of technologies in each cluster (PFS capacity sensitivity analysis, S2), hydrogen storage is measured in tonne.}
    \label{PFS_tech_s2}
% \end{figure}
% \begin{figure}[!]
\centering
    \includegraphics{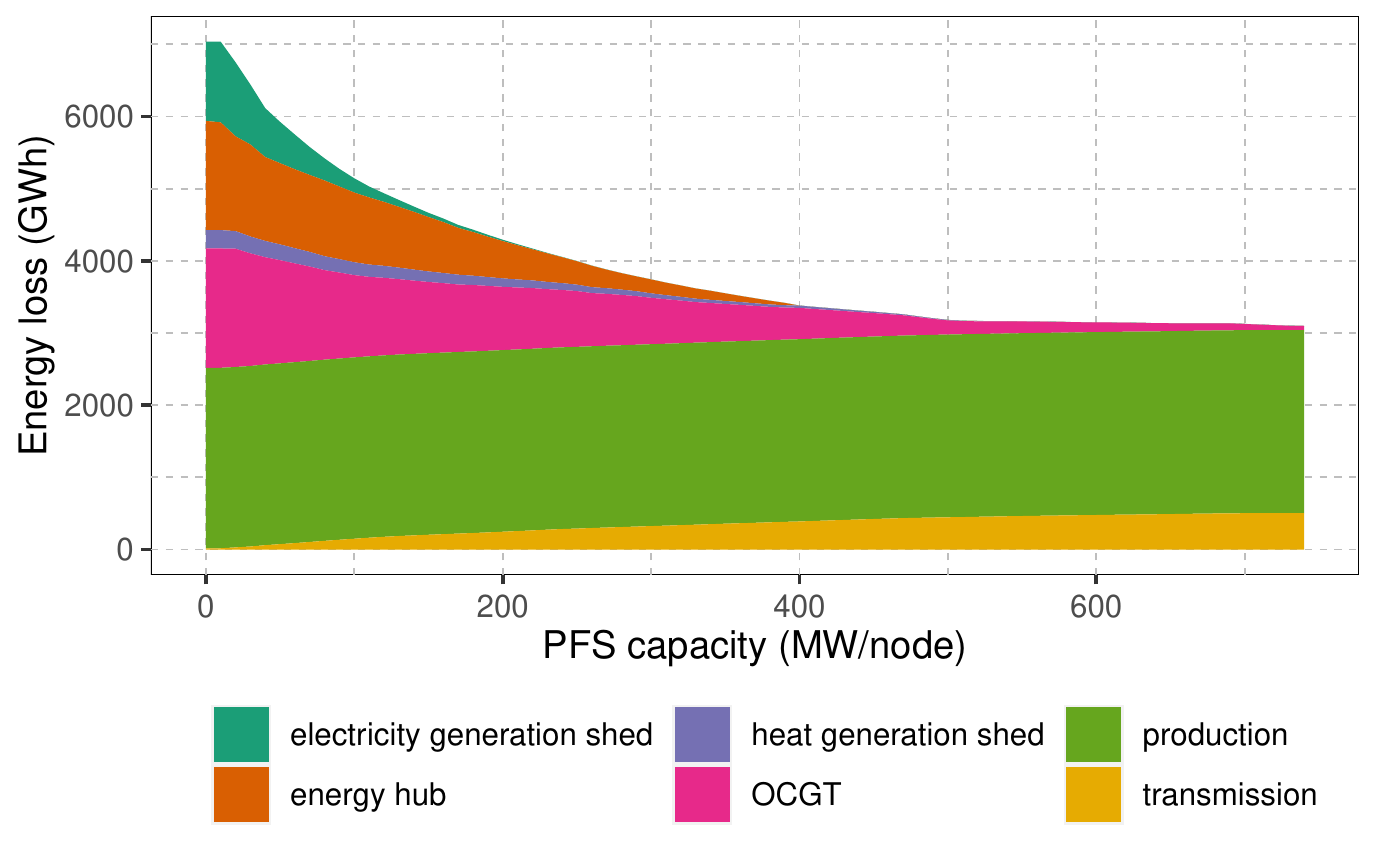}
    \caption{Energy loss (PFS capacity sensitivity analysis, S2).}
    \label{PFS_energy_loss_s2}
\end{figure}
% \begin{table}[t]
%     \caption{Combination of technologies under different AAC of PFS (S2).}
%     \centering
%     % \resizebox{\textwidth}{!}{
%     \begin{adjustbox}{width=\columnwidth,center}
% \begin{tabular}{lllll}
% \hline
% cluster  & \begin{tabular}[c]{@{}l@{}}wind\\ solar\\ energy hub\\ PFS\\ electric boiler\end{tabular} & \begin{tabular}[c]{@{}l@{}}wind\\ energy hub\\ PFS\\ electric boiler\end{tabular} & \begin{tabular}[c]{@{}l@{}}wind\\ electric boiler\\ PFS\end{tabular} & \begin{tabular}[c]{@{}l@{}}PFS\\ electric boiler\end{tabular} \\ \hline
% cluster1 &428 -- 800 &252 -- 424 &80 -- 248 &4 -- 76\\
% cluster2 &420 -- 800 & 260 -- 416 &76 -- 256 &4 -- 72\\
% cluster3 &508 -- 800, 496 -- 500 &504, 264 -- 492 &84 -- 260 &4 -- 80\\
% cluster4 &-- &292 -- 800 &76 -- 288 &4 -- 72\\
% cluster5 &360 -- 800 &276 -- 356 &84 -- 272 &4 -- 80\\ \hline
% \end{tabular}
% % }
% \label{abatement_krt_tech_table_s2}
% \end{adjustbox}
% \end{table}

In the second scenario, we fix the CO$_2$ tax to $400$ €/tonne. We increase the onshore capacity from 0 MW to 1000 MW and present the results until 740 MW.
From \ref{PFS_emission_cost_s2}, we can see that without PFS, the system can achieve 0.63 Mt/yr emissions under S2 condition. Increasing the onshore capacity brings down 43\% of the cost and also cut emission further to near zero. Figure \ref{PFS_energy_loss_s2} shows that about 20\% of the energy loss is from OEHs initially. OEHs are not needed when the onshore capacity increase to around 400 MW for each node. By adding the installed PFS capacity shown in Figure \ref{PFS_tech_s2}, we find that a total onshore capacity of 1.75 GW may help the offshore energy system achieve near zero emission. We notice that the onshore system needs to provide an averagely of 1.4 GW By checking the average power transmission of PFS, which might not be feasible without onshore system expansion.
% From Figure \ref{aac_tech_cap_s2}, we find that the AAC needs to reduce to below $252$ €/tonne to be more cost-efficient than OEHs. Figure \ref{aac_energy_loss_s2} shows that energy loss is at a minimum of around $2807$ GWh when AAC is $84$ €/tonne, and the energy is supplied by offshore wind and PFS. However, as AAC reduces further and becomes the only clean energy source, the energy loss increases $9\%$ compared with the lowest value. Under S2 condition, PFS can replace wind as the only power source if the AAC is below $72$ €/tonne.
% In the following, we present the results using metric, €/tonne, the results using €/MWh are presented in \ref{results_krmwh}.
% We see that as ACC drop from 800 €/tonne to 4 €/tonne, CO$_2$ marginal price firstly decreases then increases and finally decrease alongside the decrease of ACC of PFS. The increase of CO$_2$ marginal cost from 724 €/tonne to 532 €/tonne shows OEH' value. In this segment, the energy hub is gradually replaced by PFS. We can see that this replacement leads to a slight decrease in energy loss simply because the energy loss during transmission is not considered. To fully replace the energy hub and become the wind power compensation, its ACC must be within 288 €/tonne -- 252 €/tonne depending on specific clusters.
% \clearpage
\section{Conclusions and future work}
\label{conclusions}
This paper presents a multi-carrier offshore energy system investment planning optimisation model with a high degree of operational detail to find cost-optimal solutions for decarbonising NCS energy supply. The major novelties and contributions are (1) we formulate OEHs in a MILP model for large-scale offshore energy system planning; (2) we model the device-level energy consumption of the offshore platforms with hourly time resolution on a large scale; and (3) We use an integrated MILP investment and operational model to conduct a large-scale analysis of the value of OEHs in the North Sea offshore energy system transition towards decarbonised energy supply. Results from our case study indicate that (1) OEHs can reduce up to $10\%$ of the energy loss and $50\%$ of the emissions with CO$_2$ tax above 290 €/tonne; (2) OEHs can reduce energy loss by $65\%$ in a near zero emission system; (3) a carbon budget may be necessary to enable a zero emission energy system in addition to CO$_2$ tax; (4) the system cuts about 63\% of the initial emissions when CO$_2$ tax is 55 €/tonne, and approximately 84\% of the CO$_2$ emissions can be cut if CO$_2$ tax is increased to Norway's target of 200 €/tonne; and (5) it may be necessary to consider the whole NCS when doing system planning towards zero emission rather than consider each region separately and reach sub-optimality. 

Although the deterministic MILP model in this paper has led to many insights, there are possibilities for extensions. A deterministic optimisation model is not capable of representing load and supply uncertainties. Therefore, we aim to develop a stochastic optimisation model \cite{birge2011} and incorporate long-term and short-term uncertainties in future work. In addition, multiple investment stages are needed to represent the investment planning problem more realistically. Besides, we only consider using OEHs for fields decarbonisation, which makes OEHs seem less attractive than other technologies due to their high costs. However, OEHs can have various utilities such as energy provision to offshore fish farming, maritime and transport, and using the infrastructure past the lifetime of the oil and gas fields for purposes such as exporting hydrogen. These may motivate the investments in OEHs, and we aim to include some of the aspects in future. Finally, more work can be done on offshore network topology and the representation of the onshore power system.
\section*{CRediT author statement}
\textbf{Hongyu Zhang:} Conceptualisation, Methodology, Software, Validation, Formal analysis, Investigation, Visualisation, Data curation, Writing - original draft, Writing - review \& editing. \textbf{Asgeir Tomasgard:} Conceptualisation, Supervision, Writing - review \& editing, Funding acquisition. \textbf{Brage Rugstad Knudsen:} Conceptualisation, Supervision, Writing - review \& editing. \textbf{Harald G. Svendsen:} Conceptualisation, Writing - review \& editing. \textbf{Steffen J. Bakker:} Conceptualisation, Writing - review \& editing. \textbf{Ignacio E. Grossmann:} Conceptualisation, Supervision, Writing - review \& editing.
\section*{Declaration of competing interest}
The authors declare that they have no known competing financial interests or personal relationships that could have appeared to influence the work reported in this paper. 
\section*{Acknowledgements}
This work was supported by the Research Council of Norway through PETROSENTER LowEmission (project code 296207).
\appendix
\section{}
\label{nomenclature}
% \begin{table*}[!t]
  % \begin{framed}
    \printnomenclature[0.9in]
  % \end{framed}
% \end{table*}
\section{Complete operational constraints} 
\label{full_operational_model}
\begin{subequations}
\begin{alignat}{3}
%   &\text{s.t.}\quad&& \mathrlap{ f(x)= 
%     \sum_{i_{0}\in \mathcal{I}_{0}}(\sum_{p \in \mathcal{P}} \sum_{z \in \mathcal{Z}}C^{inv}_{pi_{0}}x^{inst}_{pzi_{0}}+\sum_{l \in \mathcal{L}}C^{inv}_{li_{0}}X^{lgth}_{l}x^{inst}_{li_{0}}) } \quad& \nonumber \\
%     &\quad&& \mathrlap{ \phantom{f(x) = } +\kappa \sum_{i \in \mathcal{I}}(\sum_{p \in \mathcal{P}}\sum_{z \in \mathcal{Z}}C^{fix}_{pi}x^{acc}_{pzi}+\sum_{l \in \mathcal{L}} C^{fix}_{li}X^{lgth}_{l}x^{acc}_{li}) } \quad& \label{investment_cost}\\
    % &\mathrlap{g(x,c)=\sum_{s \in \mathcal{S}^{T}}\sum_{t \in \mathcal{T}}\sum_{z \in \mathcal{Z}}W_{s}\left(\sum_{p \in \mathcal{P}}C^{P}_{p}y^{P}_{pzt}+C^{Shed}y^{LShed}_{zt}+ \right. } \notag \\ 
    % & \mathrlap{\hspace{9cm} \left. \sum_{g \in \mathcal{G}}C^{\text{CO}_2}\rho^{E}_{g}y^{G}_{gzt} \right)}& \label{SP_objective}\\
 &\mathrlap{g(x,c)=\sum_{s \in \mathcal{S}^{T}}W^{S}_{s}\left(\sum_{t \in
 \mathcal{T}}H_t\left(\sum_{z\in\mathcal{Z}^{P}}\left(\sum_{g \in \mathcal{G}}\left(C_g^{G}+\frac{C_g^{Fuel}+C^{\text{CO}_2}E_{g}^{Fuel}}{\eta_g^G}\right)p_{gzt}^{G}+\right.\right.\right.} &\notag\\
    &\quad&& \mathrlap{ \phantom{g(x_{i}}\left.\left. \left.\phantom{\frac{C_g^{Fuel}}{\eta_g^G}}C^{LShed}p_{zt}^{LShed}+C^{HLShed}p^{HLShed}_{zt}\right)
    +\sum_{z \in \mathcal{Z}^{O}}C^{ZO}_{zt}p^{ZO}_{zt}\right)\right)}&\\
    %  &\quad&& \phantom{g(x_{i},c_{i}) = }&\label{operational_obj}\\
    &\quad&& 0 \leq p_{gzt}^{G}+p_{gzt}^{ResG}\leq P_{gz}^{accG}, & g \in \mathcal{G}, z \in \mathcal{Z}^{P}, t \in \mathcal{T} \label{turbine_cap}\\
    &\quad&& -P_{gz}^{accG}G^{GR}_{g} \leq p_{gzt}^{G}+p_{gzt}^{ResG} &\notag\\
    &\quad&&\phantom{abc} -p_{gz(t-1)}^{G}-p_{gz(t-1)}^{ResG} \leq P_{gz}^{accG}G^{GR}_{g}, &\phantom{ab} g \in \mathcal{G}, z \in \mathcal{Z}^{P}, t \in \mathcal{T} \label{turbine_ramp}\\
    &\quad&& p_{rzt}^{R}=R^{R}_{rzt}p_{rz}^{accR}, &  r \in \mathcal{R}, z \in \mathcal{Z}^{H}, t \in \mathcal{T} \label{renewable_cap}\\
    &\quad&& 0 \leq p_{zt}^{ZO}\leq p_{z}^{accZO}, &  z \in \mathcal{Z}^{O},t \in \mathcal{T} \label{onshore_cap}\\
    &\quad&& q_{sz(t+1)}^{SE}=q_{szt}^{SE}+H_{t}(\eta_{s}^{SE}p_{szt}^{SE+}-p_{szt}^{SE-}), & s \in \mathcal{S}^{E}, z \in \mathcal{Z}^{P}, t \in \mathcal{T} \label{estore_balance}\\
    &\quad&& H_{t}(p_{szt}^{ResSE} + p_{szt}^{SE-})\leq q_{szt}^{SE}, & s \in \mathcal{S}^{E}, z \in \mathcal{Z}^{P}, t \in \mathcal{T} \label{estore_discharge_cap_1}\\
    &\quad&&0 \leq p_{szt}^{SE+} \leq H^{SE}_{s}q_{sz}^{accSE}, & s \in \mathcal{S}^{E}, z \in \mathcal{Z}^{P}, t \in \mathcal{T} \label{estore_charge_cap}\\
    &\quad&&0 \leq p_{szt}^{SE-} + p_{szt}^{ResSE} \leq H^{SE}_{s}q_{sz}^{accSE}, & s \in \mathcal{S}^{E}, z \in \mathcal{Z}^{P}, t \in \mathcal{T} \label{estore_discharge_cap}\\
    &\quad&&0 \leq q_{szt}^{SE} \leq q_{sz}^{accSE}, & s \in \mathcal{S}^{E}, z \in \mathcal{Z}^{P}, t \in \mathcal{T} \label{estore_energy_cap}\\
    &\quad&&\sum_{g \in \mathcal{G}}p_{gzt}^{G}+
     \sum_{l \in \mathcal{L}}A^{E}_{zl}\eta^{L}_{l}p_{lt}^{L}& \notag\\
     &\quad&&\phantom{abc}+\sum_{s \in \mathcal{S}^{E}}p_{szt}^{SE-}+p_{zt}^{LShed}=& \notag\\
    &\quad&&\phantom{abc}\mu^{D}p_{zt}^{D}+\sum_{s \in \mathcal{S^{E}}}p_{szt}^{SE+}+
     p_{zt}^{GShed}, &  z \in \mathcal{Z}^{P}, t \in \mathcal{T}\label{kcl_platform}\\
    &\quad&&\sum_{r \in \mathcal{R}}p_{rzt}^{R}+
     \sum_{l \in \mathcal{L}}A^{E}_{zl}\eta^{L}_{l}p_{lt}^{L}+\sum_{f \in \mathcal{F}}p^{F}_{fzt}=&\notag\\
    &\quad&&\phantom{abc}p_{zt}^{GShed}+\sum_{e \in \mathcal{E}}p^{E}_{ezt},& z \in \mathcal{Z}^{H}, t \in \mathcal{T}\label{kcl_energyhub}\\
    &\quad&&p_{zt}^{ZO}+\sum_{l \in \mathcal{L}}A^{E}_{zl}\eta^{L}_{l}p_{lt}^{L}= p_{zt}^{GShed},
        & z \in \mathcal{Z}^{O}, t \in \mathcal{T}\label{kcl_onshore}\\
    &\quad&&-p_{l}^{accL} \leq p_{lt}^{L} \leq p_{l}^{accL}, & l \in \mathcal{L}, t \in \mathcal{T}\label{line_cap}\\
    &\quad&&\sigma^{Res}p^{D}_{zt} \leq \sum_{g \in \mathcal{G}} p^{ResG}_{gzt}+ \sum_{s \in \mathcal{S}^{E}}p^{ResSE}_{szt}, & z \in \mathcal{Z}^{P}, t \in \mathcal{T}\label{reserve}\\
    &\quad&&p^{D}_{zt}=\sum_{b \in \mathcal{B}^{E}}p^{EB}_{bzt}+\sum_{c \in \mathcal{C}^{Exp}}p^{CExp}_{czt} &\notag\\
    &\quad&&\phantom{abc}+\sum_{c \in \mathcal{C}^{Inj}}p^{CInj}_{czt}+\sum_{p \in \mathcal{P}^{O}}p^{PO}_{pzt} & \notag\\
    &\quad&&\phantom{abc}+\sum_{p \in \mathcal{P}^{WI}}p^{PWI}_{pzt}+\sum_{p \in \mathcal{P}^{WL}}p^{PWL}_{pzt}, & z \in \mathcal{Z}^{P}, t \in \mathcal{T}\label{power_demand}\\
    &\quad&&\mathrlap{\sum_{s \in \mathcal{S}^{T}}W^{S}_{s}\sum_{z \in \mathcal{Z}^{P}} \sum_{t \in \mathcal{T}} \sum_{g \in \mathcal{G}}\frac{E_{g}^{Fuel}p_{gzt}^{G}H_t}{\eta_{g}^{G}} \leq \mu^{E},}\label{co2budget}&
    \end{alignat}
    \end{subequations}
     \begin{subequations}
     \label{production_constraints}
    \begin{alignat}{3}
    &\phantom{\text{sss.}}\quad&&0 \leq p^{EB}_{bzt} \leq p_{bz}^{accEB},&  b \in \mathcal{B}^{E}, z \in \mathcal{Z}^{P},     t \in \mathcal{T}\label{eboiler_cap}\\
    &\quad&&0 \leq p^{HSEP}_{dzt} \leq P^{accSEP}_{dz}, & d \in \mathcal{D}^{Sep}, z \in \mathcal{Z}^{P}, t \in \mathcal{T} \label{separator_cap}\\
    &\quad&&\sum_{g \in \mathcal{G}}\eta^{Gh}_{g}p^{G}_{gzt}+\sum_{b\in \mathcal{B}^{E}}\eta^{Eb}_{b}p^{EB}_{bzt}+p^{HLShed}_{zt} & \notag\\
    &\quad&&\phantom{abc}=\sum_{d \in \mathcal{D}^{Sep}}p^{HSEP}_{dzt}+p^{HGShed}_{zt},& z \in \mathcal{Z}^{P}, t \in \mathcal{T}\label{heat_balance}\\
    &\quad&&p^{HSEP}_{dzt}=\rho^{SEP}_{d}V^{OD}_{zt}, & d \in \mathcal{D}^{Sep},  z \in \mathcal{Z}^{P}, t \in \mathcal{T}\label{separator_heat}\\
    &\quad&&0 \leq p_{czt}^{CExp} \leq P_{cz}^{accCExp}, & c \in \mathcal{C}^{Exp}, z \in \mathcal{Z}^{P}, t \in \mathcal{T}\label{compressor_cap1}\\
    &\quad&& 0 \leq p_{czt}^{CInj} \leq P_{cz}^{accCInj}, & c \in \mathcal{C}^{Inj}, z \in \mathcal{Z}^{P}, t \in \mathcal{T}\label{compressor_cap2}\\
    &\quad&& p_{czt}^{CInj}=\frac{V_{zt}^{NGI}}{\eta^{C}H\rho}[(\gamma_{c})^{\frac{\alpha-1}{\alpha}}-1], \quad& c \in \mathcal{C}^{Inj}, z \in \mathcal{Z}^{P}, t \in \mathcal{T}\label{compressor_power1}\\
    &\quad&& p_{czt}^{CExp}=\frac{V_{zt}^{NGE}}{\eta^{C}H\rho}[(\gamma_{c})^{\frac{\alpha-1}{\alpha}}-1], \quad& c \in \mathcal{C}^{Exp}, z \in \mathcal{Z}^{P}, t \in \mathcal{T}\label{compressor_power2}\\
    &\quad&& 0 \leq p_{pzt}^{PO} \leq P_{pz}^{accPO}, &  p \in \mathcal{P}^{O}, z \in \mathcal{Z}^{P}, t \in \mathcal{T}\label{pump_oil_cap}\\
    &\quad&&p_{pzt}^{PO} = \kappa^{PO}_{p}V^{OD}_{zt}, &  p \in \mathcal{P}^{O}, z \in \mathcal{Z}^{P}, t \in \mathcal{T}\label{pump_oil_power}\\
    &\quad&&0 \leq p_{pzt}^{PWI} \leq P_{pz}^{accPWI}, & p \in \mathcal{P}^{WI}, z \in \mathcal{Z}^{P}, t \in \mathcal{T}\label{pump_water_inj_cap}\\
    &\quad&&0 \leq p_{pzt}^{PWL} \leq P_{pz}^{accPWL},& p \in \mathcal{P}^{WL}, z \in \mathcal{Z}^{P}, t \in \mathcal{T}\label{pump_water_lift_cap}\\
    &\quad&&p_{pzt}^{PWI} = \kappa^{PWI}_{p}V^{WI}_{zt}, & p \in \mathcal{P}^{WI}, z \in \mathcal{Z}^{P}, t \in \mathcal{T}\label{pump_water_inj_power}\\
    &\quad&&p_{pzt}^{PWL} = \kappa^{PWL}_{p}V^{WL}_{zt}, & p \in \mathcal{P}^{WL}, z \in \mathcal{Z}^{P}, t \in \mathcal{T}\label{pump_water_lift_power}
    \end{alignat}
    \end{subequations}
     \begin{subequations}
     \label{OEH_constraints}
    \begin{alignat}{3}
    &\phantom{\min}\quad&&0\leq p_{ezt}^{E} \leq p_{ez}^{accE}, &\qquad e \in \mathcal{E}, z \in \mathcal{Z}^{H}, t \in \mathcal{T}\label{electrolyser_cap}\\
    &\quad&&0 \leq p_{fzt}^{F} \leq p_{fz}^{accF},& f \in \mathcal{F}, z \in \mathcal{Z}^{H},  t \in \mathcal{T}\label{fuel_cell_cap}\\
    &\quad&&0 \leq v_{szt}^{SHy} \leq v_{sz}^{accSHy},&   s \in \mathcal{S}^{Hy}, z \in \mathcal{Z}^{H},  t \in \mathcal{T}\label{hydrogen_cap}\\
    &\quad&&|p_{fzt}^{F}- p_{fz(t-1)}^{F}|\leq F^{FR}_{f}p_{fz}^{accF}, &   f \in \mathcal{F}, z \in \mathcal{Z}^{H},  t \in \mathcal{T}\label{fuel_cell_ramp}\\
    &\quad&&v^{SHy}_{sz(t+1)}=v^{SHy}_{szt}+v^{SHy+}_{szt}-v^{SHy-}_{szt}, & s \in \mathcal{S}^{Hy}, z \in \mathcal{Z}^{H}, t \in \mathcal{T}\label{hydrogen_storage_balance}\\
    &\quad&&
    \eta^{EF}(\sum_{f \in \mathcal{F}}\frac{p^{F}_{fzt}H_{t}}{\eta^{F}_{f}\theta^{Hy}}-\sum_{s \in \mathcal{S}^{Hy}}v^{SHy-}_{szt})=&\notag\\
    &\quad&&\phantom{abc}\sum_{e \in \mathcal{E}}p^{E}_{ezt}H_{t}-\eta^{ES}\sum_{s \in \mathcal{S}^{Hy}}v^{SHy+}_{szt}, &   z \in \mathcal{Z}^{H},  t \in \mathcal{T}.\label{hydrogen_balance}
\end{alignat}
\end{subequations}
\eqref{turbine_cap}--\eqref{co2budget} give the electricity system constraints. \eqref{eboiler_cap}--\eqref{pump_water_inj_power} give the platform production constraints. \eqref{electrolyser_cap}--\eqref{hydrogen_balance} represent hydrogen system constraints of OEHs. All variables except for power transmission is non-negative. Power transmission is a free variable. All variables are also indexed by operational year $i$, and we omit them for ease of notation.
\section{Calculation of energy loss}
The indices, summation and multiplication of one hour are omitted.
\label{sec:energy_loss}
\begin{equation}
\begin{split}
    q^{Loss}=&p^{GShed}+p^{HGShed}+(\frac{1}{\eta^{G}}-1-\eta^{Gh})p^{G}+(1-\eta^{SE})p^{SE+}\\
    &+(1-\eta^{SEP})p^{HSEP}+(1-\eta^{C})(p^{CExp}+p^{CInj})+(1-\eta^{PO})p^{O}\\
    &+(1-\eta^{PWL})p^{PWL}+p^{E}-\theta^{Hy}(\frac{p^{F}}{\eta^{F}\theta^{Hy}}-v^{SHy-}+v^{SHy+})\\
    &+(\frac{1}{\eta^{F}}-1)p^{F}+(1-\eta^{L})p^{L}+(1-\eta^{PWI})p^{PWI}+(1-\eta^{EB})p^{EB}.
\end{split}\notag
\end{equation}
% \section{Visualisation of the results}
% \begin{figure}[!]
%     \centering
%     \includegraphics{abatement_krt_energy_loss_s2.pdf}
%     \caption{Energy loss (sensitivity analysis of AAC, S2).}
%     \label{aac_energy_loss_s2}
% \end{figure}
% \begin{figure}[!]
%     \centering
%     \includegraphics{abatement_krt_cost_s2.pdf}
%     \caption{Emission and cost (sensitivity analysis of AAC, S2).}
%     \label{abatement_krt_emission_energyloss_s1}
% \end{figure}
% \section{Results of sensitivity analysis of ACC of PFS using kr/MWh}
% \label{results_krmwh}
% \begin{figure}[]
%     \centering
%     \caption{Emission and energy loss (sensitivity analysis of ACC (kr/MWh), S1).}
%     \includegraphics{abatement_krmwh_emission_energyloss_s1.pdf}
%     \label{abatement_krmwh_emission_energyloss_s1}
% \end{figure}
% \begin{figure}[]
%     \centering
%     \caption{Capacities of technologies in each cluster (sensitivity analysis of ACC (kr/MWh), S1).}
%     \includegraphics{abatement_krmwh_tech_cap_s1.pdf}
%     \label{abatement_krmwh_tech_cap_s1}
% \end{figure}
% \begin{figure}[]
%     \centering
%     \caption{Energy loss and CO$_2$ marginal price (sensitivity analysis of ACC (kr/MWh), S2).}
%     \includegraphics{abatement_krmwh_emission_energyloss_s2.pdf}
%     \label{abatement_krmwh_emission_energyloss_s2}
% \end{figure}
% \begin{figure}[]
%     \centering
%     \caption{Capacities of technologies in each cluster (sensitivity analysis of ACC (kr/MWh), S2).}
%     \includegraphics{abatement_krmwh_tech_cap_s2.pdf}
%     \label{abatement_krmwh_tech_cap_s2}
% \end{figure}
\clearpage





% \section*{References}
%\bibliographystyle{elsarticle-harv}
\setlength{\bibsep}{0pt plus 0.3ex}
\footnotesize{
\bibliographystyle{model5-names}%When changing reference style
\bibliography{literature}}
\end{document}